\documentclass[12pt]{article}
\UseRawInputEncoding
\usepackage{amssymb}
\usepackage{latexsym}
\usepackage{amsfonts}
\usepackage{amsmath}
\newtheorem{theorem}{Theorem}[section]
\newtheorem{lemma}[theorem]{Lemma}
\newtheorem{proposition}[theorem]{Proposition}

\newtheorem{corollary}[theorem]{Corollary}
\newtheorem{remark}[theorem]{Remark}

\numberwithin{equation}{section}

\usepackage[dvips]{color}

\def\P{\mathbb{P} }
\def\E{\mathbb{E} }
\def\R{\mathbb{R} }

\def\tl{\tilde}
\def\ocM{\overline{{\cal M}}_{r,t}}
\def\ucM{\underline{{\cal M}}_{r,t}}
\def\ocL{\overline{{\cal L}}_{r,t}}
\def\ucL{\underline{{\cal L}}_{r,t}}

\def\mP{{\mbox{P}} }
\def\mE{{\mbox{E}} }

\def\cL{{\cal L} }
\def\cE{{\cal E} }

\def\cM{{\cal M}}
\def\cS{{\cal S}}

\begin{document}
\allowdisplaybreaks

\title{\Large\bf   The extremal process of super-Brownian motion}
\author{ \bf  Yan-Xia Ren\footnote{The research of this author is supported by NSFC (Grant No. 11671017  and 11731009) and LMEQF.\hspace{1mm} } \hspace{1mm}\hspace{1mm}
Renming Song\thanks{Research supported in part by a grant from the Simons
Foundation (\#429343, Renming Song).} \hspace{1mm}\hspace{1mm} and \hspace{1mm}\hspace{1mm}
Rui Zhang\footnote{
The corresponding author.
The research of this author is supported by NSFC (Grant No. 11601354), Beijing Municipal Natural Science Foundation(Grant No. 1202004), and Academy for Multidisciplinary Studies, Capital Normal University}
\hspace{1mm} }
\date{}
\maketitle

\begin{abstract}
In this paper, we establish limit theorems for the
supremum of the support, denoted by $M_t$,
of a supercritical super-Brownian motion $\{X_t, t\ge0\}$ on $\R$.
We prove that there exists an $m(t)$ such that $(X_t-m(t), M_t-m(t))$ converges in law,
and give some large deviation results for $M_t$ as $t\to\infty$.
We also prove that
the limit of the extremal process $\mathcal{E}_t:=X_t-m(t)$ is a
Poisson random measure with exponential intensity in which each atom
is decorated by an independent copy of an auxiliary measure.
These results are analogues of the results for branching Brownian motions
obtained in Arguin et al. (Probab. Theory Relat. Fields \textbf{157} (2013), 535--574), A\"id\'ekon et al.
(Probab. Theory Relat. Fields \textbf{157} (2013), 405--451) and Roberts (Ann. Probab. \textbf{41} (2013),  3518--3541).
\end{abstract}

\medskip

\noindent{\bf Keywords and Phrases}: Super-Brownian motion, extremal process,
supremum of the support of super-Brownian motion, Poisson random measure,  KPP equation.

\begin{doublespace}
\section{Introduction}

\subsection{Super-Brownian motion}\label{ss:SBM}

Let $\psi$ be a function of the form:
$$\psi(\lambda)=-\alpha\lambda+\beta\lambda^2+\int_0^\infty \Big(e^{-\lambda y}-1+\lambda y\Big)n(dy),\quad \lambda\ge 0,$$
where
$\alpha\in\R$,
$\beta \ge 0$ and $n$ is a $\sigma$-finite measure satisfying
$$\int_0^\infty (y^2\wedge y) n(dy)<\infty.$$
$\psi$ is called a branching mechanism.
We will always assume that $\lim_{\lambda\to\infty}\psi(\lambda)=\infty$.
Let $\{B_t,\mbox{P}_x\}$ be a standard Brownian motion, and $\mE_x$ be the corresponding expectation.
In this paper we will consider a super-Brownian motion $X$ on $\R$ with branching mechanism $\psi$.

Let $\mathcal{B}^+(\R)$ (resp. $\mathcal{B}^+_b(\R)$) be the space of non-negative
(resp. non-negative bounded) Borel measurable function on $\R$,
and let ${\cal M}_F(\R)$
be the space of finite measures on $\R$, equipped with the topology of weak convergence.
A super-Brownian motion $X$ with branching mechanism $\psi$ is a  Markov process taking values in ${\cal M}_F(\R)$.
The existence of such superprocesses is well-known, see, for instance,
\cite{Dawson}, \cite{E.B.} or \cite{Li11}. For any $\mu \in \mathcal{M}_F(\R)$, we denote the
law of $X$ with initial configuration $\mu$
 by $\P_\mu$, and the corresponding expectation by $\E_\mu$.
As usual, we use the notation:
$\langle \phi,\mu\rangle:=\int_{\R} \phi(x)\mu(dx)$
and $\|\mu\|:=\langle 1,\mu\rangle$. Then for all
$\phi\in \mathcal{B}^+_b(\R)$ and $\mu \in \mathcal{M}_F(E)$,
\begin{equation}
    -\log \E_\mu\left(e^{-\langle \phi,X_t\rangle}\right)=
  \langle u_\phi(t, \cdot),\mu\rangle,
\end{equation}
where $u_\phi(t, x)$ is the unique positive solution to the equation
\begin{equation}\label{eqt-u}
  u_\phi(t, x)+\mE_x\int_0^t\psi(u_\phi(t-s, B_s))ds=\mE_x \phi(B_t).
\end{equation}
Note that the integral equation \eqref{eqt-u} is equivalent to  the equation:
\begin{equation}\label{KPPpsi}
  \frac{\partial}{\partial t}u_{\phi}(t,x)-\frac{1}{2}\frac{\partial^2}{\partial x^2}u_{\phi}(t,x)=-\psi(u_{\phi}(t,x)),\quad t>0, x\in\R,
\end{equation}
with initial condition $u_{\phi}(0,x)=\phi(x)$.
Moreover, $\lim_{t\to0}u_\phi(t,x)=\phi(x)$, if
$\phi$ is a nonnegative bounded continuous function on $\R$.

$X$ is called a supercritical (critical, subcritical) super Brownian motion if $\alpha>0\, (=0,<0)$. In this paper, we only deal with the supercritical case, that is $\alpha>0$.

\subsection{Maximal position of super-Brownian motion}\label{ss:maxofSBM}

The maximal  position $M_t$ of branching-Brownian motions
has been studied intensively.
Without loss of generality, we assume in this subsection that the branching rate is $1$,
and the offspring distribution $\{p_k\}$ satisfies $p_0=0$ and the mean of the offspring distribution is $2$.
Denote by ${\bf P}_{\delta_0}$ the law of branching Brownian motion
starting from one point located at $0$.
In the seminal paper \cite{KPP}, Kolmogorov,  Petrovskii and  Piskounov proved that
$M_t/t\to \sqrt{2}$ in probability, which implies that the leading order of $M_t$ is $\sqrt{2} t.$
In \cite{Bramson78}, Bramson provided a log correction to the leading order  of $M_t$.
He proved in \cite{Bramson78} (see also \cite{Bramson}) that, under some moment conditions,
${\bf P}_{\delta_0}(M_t-m(t)\le x)\to 1-w(x)$
as $t\to\infty$ for all $x\in \R$, where $m(t)=\sqrt{2}t-\frac{3}{2\sqrt{2}}\log t$ and  $w(x)$ is a traveling wave solution.
In \cite{LS87}, Lalley and Sellke gave
a probabilistic representation of the traveling wave solution
in terms of the limit of the derivative martingale of branching Brownian motion.
In \cite{Robert}, Roberts gave another proof of Bramson's result and also an almost
sure fluctuation result of $M_t$.
Large deviation results for $M_t$ were obtained by Chauvin and Rouault in \cite{Chauvin88,Chauvin}.

Beyond the behavior of the maximal displacement of branching Brownian motions,
the full statistics of the extremal configurations was studied in Arguin et al. \cite{ABK11, ABK12,ABK} and A\"{i}d\'ekon et al. \cite{ABBS}.
Assume the particles alive at time $t$ are ordered decreasingly: $x^t_{1}\ge x^t_{2}\ge \cdots\ge x^t_{n(t)}$, where $n(t)$ is the number of particles alive at time $t$. It is clear that $x_1^t$ is the maximum position $M_t$ at time $t$.
Arguin et al. \cite{ABK12,ABK} studied the limit property of the extremal process of branching Brownian motion, which is the random measure defined by
$$
\cE_t:=\sum_{j=1}^{n(t)}\delta_{x^t_j-m(t)}.
$$
Note that $\cE_t=Y_t-m(t)$, where $Y_t$ is the measure corresponding to configuration
of the positions of the particles alive at time $t$.
In \cite{ABK},
using the results of \cite{Bramson},  Arguin et al. first proved that $\cE_t$ converges in law,
which implies  the weak convergence of $x^t_{k}$, the $k$th  maximal displacement for each fixed integer $k\geq 1$, and then
gave a rigorous  characterization of the limiting extremal process.
It was proved in  \cite{ABK} that the limiting
process is a (randomly shifted) Poisson cluster process, where the positions of the
clusters form a Poisson point process with an exponential intensity measure.
The law of the individual clusters is characterized as
a branching Brownian motion conditioned
to perform unusually large displacements.  Almost at the same time,
A\"{i}d\'ekon et al. \cite{ABBS} proved similar results using a totally different method.

In the recent paper \cite{BBCM}, Berestycki et al.  studied the asymptotic behavior of the extremal particles of branching Ornstein-Uhlenbeck processes.
 For  inhomogeneous branching Brownian motions, many papers discussed the growth rate of the maximal position, see Bocharov and Harris \cite{Bocharov-Harris14,Bocharov-Harris16} and Bocharov \cite{Bocharov} for the case with catalytic branching at the origin, Shiozawa \cite{Shiozawa}, Nishimori et al. \cite{Nishimori-Shiozawa}, Lalley and Sellke \cite{LS88,LS89} for the case with some general branching mechanisms.
  For  branching random walks,
 we refer the readers to Hu et al. \cite{HS}, A\"{i}d\'{e}kon \cite{Aldekon}, Madaule \cite{Madaule} and Carmona et al. \cite{Carmona-Hu}.

Unlike the case of branching Brownian motions or branching random walks, there are very few
results for the supremum of super-Brownian motions, see \cite{KLMR, Englander}.
Let $X_t$ be the super-Brownian motion in Subsection \ref{ss:SBM} and let $M_t$ be the supremum of the support of $X_t$. We will prove that, under some
conditions, $\P_{\delta_0}(M_t-m(t)\le x)\to e^{-w(x)}$ as $t\to\infty$ for all $x\in \R$,
where $m(t):=\sqrt{2\alpha} t-\frac{3}{2\sqrt{2\alpha}}\log t$ and
$w$ is a traveling  wave  solution.
We also give some large deviation results for $M_t$.
In analogy to the case of branching Brownian motions, we will call the random measure $\cE_t:=X_t-m(t)$ the extremal process of the super-Brownian motion
$X$, which is simply the super-Brownian motion seen from the position $m(t)$.
We will generalize the results in \cite{ABK} to super-Brownian motions and study the limit
of $\cE_t$.
We will give the precise statements
of our main results in Subsection \ref{ss:mainresults}.

Our proofs depend heavily on the convergence of solutions of the Kolmogorov-Petrovsky-Piscounov (KPP) equation \eqref{KPPpsi}, with general initial conditions not necessarily bounded between 0 and 1, to traveling wave solutions.

\subsection{KPP equation related to super-Brownian motion}\label{ss:KPP}

The classical KPP equation is a semilinear equation of the form
\begin{equation}\label{KPPequ}
  u_t(t, x)-\frac{1}{2}u_{xx}(t, x)=f(u(t, x)),\quad (t,x)\in(0,\infty)\times \R.
\end{equation}
The KPP equation has been studied for many years  analytically, see for example,
Kolmogorov et al. \cite{KPP}, Fisher \cite{Fisher}, Aronson et al. \cite{Aronson}, Bramson \cite{Bramson}, Lau \cite{Lau}, Volpert et al. \cite{VV}.

In \cite{Bramson},
the nonlinear function $f$ can be any function on $[0,1]$ satisfying
\begin{equation}\label{f-condition1}
\begin{array}{l}
f\in C^1[0, 1], f(0)=f(1)=0, f(u)>0  \mbox{ for } u\in (0, 1); \\
f'(0)=1, f'(u)\le 1, \mbox{ for } 0<u\le 1,
\end{array}
\end{equation}
and
\begin{equation}\label{f-condition2}
1-f'(u)=O(u^\rho)\quad (\mbox{as }u\to 0) \mbox{ for some } \rho>0.
\end{equation}
Kolmogorov et al. \cite{KPP} showed that under condition
\eqref{f-condition1} and with Heaviside initial  condition $u(0,x)=1_{(-\infty, 0)}(x)$, \eqref{KPPequ}
 has a unique solution $u(t, x)$ satisfying
\begin{equation}\label{bramson}
  \lim_{t\to\infty}u(t,m(t)+x)=w(x),\quad \mbox{uniformly in } x\in\R,
\end{equation}
for some centering term $m(t)$, where $m(t)$ satisfies $m(t)=\sqrt{2}t+o(t)$ as $t\to\infty$, and $w$ is a travelling wave solution,
which is a function solving
the ordinary differential equation $\frac{1}{2}w_{xx}+\sqrt{2}w_x+f(w)=0$, and satisfying $0<w(x)<1$, $\lim_{x\to\infty}(x)=0, \lim_{x\to-\infty}w(x)=1$.
Bramson \cite{Bramson} improved the above result in
two aspects: first the initial
 condition $u(0,x)$ is a general function between $0$ and $1$,
not just the Heaviside initial  condition $u(0,x)=1_{(-\infty, 0)}(x)$;
 secondly he proved that if in addition $f$ satisfies \eqref{f-condition2} and  the initial condition
 $u(0,x)$ satisfy some integrability condition, \eqref{bramson} holds with $m(t)=\sqrt{2}t-\frac{3}{2\sqrt{2}}\log t$.
Note that, since $0$ and $1$ are two special solutions,
it follows from the maximum principle that
any solution of \eqref{KPPequ}, with
initial condition bounded between 0 and 1, must be bounded between 0 and 1.

An interesting link between branching Brownian motion  and partial differential equations was observed
 by McKean \cite{McKean} (see also Ikeda, Nagasawa and Watanabe \cite{INW1, INW2, INW3}):
$u(t,x):={\bf P}_{\delta_0}(M_t> x)$
solves the KPP equation \eqref{KPPequ}
with initial condition $u(0,x)=1_{(-\infty,0)}(x)$ and with
$f(u)=(1-u)-\sum^\infty_{k=0}p_k(1-u)^k$,
where $\{p_k, k\geq 0\}$ is the offspring distribution and the branching rate is 1.
Moreover, if $p_0=0$, $\sum_{k}kp_k=2$, and $\sum_{k}k^{1+\rho}p_k<\infty$, then $f(u)$ satisfies conditions \eqref{f-condition1} and \eqref{f-condition2}.
In probabilistic language, \eqref{bramson} gives the convergence in distribution for $M_t-m(t)$.
There are also some papers using branching Brownian motions to study travelling wave solutions to the KPP equation, see \cite{Harris,Kyprianou}, for instance.

It follows from \eqref{KPPpsi} that the super-Brownian motion $X$ is related to the KPP equation with $f=-\psi$. It is natural to use this relationship to investigate the maximal position of super-Brownian motions.
Let $\lambda^*$ be the largest root of the equation $\psi(\lambda)=0$. Since $\psi'(0)=-\alpha<0$, $\psi(\infty)=\infty$, it follows from the strict convexity of $\psi$ that $\lambda^*>0$ exists.
Note that $0$ and $\lambda^*$ are two special solutions of \eqref{KPPpsi}.
One might think that the role of $0$ and  $\lambda^*$ for the KPP \eqref{KPPpsi} corresponding to super-Brownian motions is similar that of $0$ and $1$ for the KPP equation \eqref{KPPequ}
corresponding to branching Brownian motions.
However, for super-Brownian motions we need to consider general non-negative solutions of the corresponding KPP equation \eqref{KPPpsi} with initial condition  $u(0,x)$ not necessarily bounded between $0$ and $\lambda^*$.
In this paper, we will first generalize Bramson's results in \cite{Bramson}  to
general non-negative solutions of the  KPP equation \eqref{KPPpsi} associated with super-Brownian motions, with initial conditions not necessarily bounded between 0 and $\lambda^*$,
see  \eqref{heavy-boundary} below for example.
Let $u_{\phi}(t, x)$ be a non-negative solution to \eqref{KPPpsi} with initial condition $\phi$.
In this paper, we will prove that there also exists some function $m(t)$ such that,
for general initial condition $\phi$, $u_{\phi}(t,m(t)+x)$ converges to some traveling wave solution. More precisely, we consider non-increasing traveling wave solutions $w$ with speed $\sqrt{2\alpha}$ to the equation \eqref{KPPpsi} such that
$$
  \lim_{x\to\infty}w(x)=0,\qquad \lim_{x\to -\infty}w(x)=\lambda^*.
$$
By a non-increasing traveling wave solution with speed $\sqrt{2\alpha}$ to  \eqref{KPPpsi},
 we mean a non-negative non-increasing function $w$ such that
 $w(x-\sqrt{2\alpha}t)$  is a solution to \eqref{KPPpsi}.
 Clearly, $w$ satisfies
 $$
  \frac{1}{2}w_{xx}+\sqrt{2\alpha}w_x-\psi(w)=0.
 $$
 We will give an exact asymptotic expression for $m(t)$.
We will then use these results to study asymptotic properties of the
supremum of the support and the extremal process of the super-Brownian motion $X$.

\subsection{Main results}\label{ss:mainresults}
We will assume that $\psi$ satisfies the following two conditions:
\begin{itemize}
  \item[] {\bf{(H1)}} There exists $\gamma>0$ such that
  \begin{equation}\label{cond-log}
   \int_1^\infty y(\log y)^{2+\gamma}n(dy)<\infty.
  \end{equation}
  \item[] {\bf{(H2)}} $\psi$ satisfies
\begin{equation}\label{e:KLMR(5)}
\int^\infty\frac1{\sqrt{\int^{\xi}_{\lambda^*}\psi(u)du}}d\xi<\infty.
\end{equation}
\end{itemize}

Let $\mathcal{R}$ be the smallest closed set such that $\mbox{supp }X_t\subseteq\mathcal{R}$, $t\ge 0$.
It is known (cf. \cite{Sheu}) that {\bf{(H2)}} implies Grey's condition
\begin{equation}\label{Grey}
  \int^\infty \frac{1}{\psi(\lambda)}d\lambda<\infty
\end{equation}
and that
$$\P_\mu(\mathcal{R} \mbox{ is compact})=e^{-\lambda^*\|\mu\|}.$$
It is well known that $\{\|X_t\|\}$ is a continuous state branching process and that, under condition $\eqref{Grey}$,
\begin{align}\label{extinction}
  \P_{\mu}(\|X_t\|=0)>0
\end{align}
and
$\lim_{t\to\infty}\P_{\mu}(\|X_t\|=0)=e^{-\lambda^*\|\mu\|}.$
Denote  $\mathcal{S}:=\{\forall t\ge 0, \|X_t\|>0\}$.

For some of our results, we also need the following stronger assumption:
\begin{itemize}
  \item[]{\bf (H3)}
There exist $\vartheta\in(0,1]$ and $a>0,b>0$ such that
$$
\psi(\lambda)\ge -a\lambda+b\lambda^{1+\vartheta}, \quad \lambda >0.
$$
\end{itemize}
Clearly, condition $\bf{(H3)}$ implies $\bf{(H2)}$. In particular, $\bf{(H3)}$ holds if $\beta>0$.
Actually, condition $\bf{(H3)}$ is only used in proving Lemma \ref{initial-cond-V}.

 Note that super-Brownian motions have been used to study traveling wave solutions to the KPP equation \eqref{KPPpsi}, see \cite{Kyprianou,KLMR}, for instance.
For convenience,
we write $\P:=\P_{\delta_0}$ and $\E:=\E_{\delta_0}$.
Define, for $t\ge0$,
$$
Z_t:=\langle (\sqrt{2\alpha}t-\cdot)e^{-\sqrt{2\alpha}(\sqrt{2\alpha}t-\cdot)},X_t \rangle.
$$
It has been proven in \cite{KLMR}
that $\{Z_t, t\geq 0\}$ is a martingale,
which is called the derivative martingale,
and that $Z_t$ has an almost sure non-negative limit $Z_\infty$ as $t\to\infty$.
Furthermore, $Z_\infty$ is almost surely  positive on
$\cS$
if and only if
\begin{equation}\label{cond:n}
  \int_1^\infty y (\log y)^2n(dy)<\infty.
\end{equation}
Clearly, \eqref{cond-log} implies \eqref{cond:n}. Thus $Z_\infty$ is almost surely  positive on $\cS$. The traveling wave solution with speed $\sqrt{2\alpha}$ to  \eqref{KPPpsi} is given by
\begin{equation}\label{Trav-Solution}
  w(x)=-\log \E\left[\exp\left\{-cZ_\infty e^{-\sqrt{2\alpha}x}\right\}\right]
\end{equation}
and
\begin{equation}\label{limit-w}
  \lim_{x\to\infty}\frac{w(x)}{xe^{-\sqrt{2\alpha}x}}=c.
\end{equation}
For more details, we refer our readers to \cite[Theorems 2.4 and 2.6]{KLMR}.
Under condition \eqref{cond:n}, $Z_\infty=\Delta(\sqrt{2\alpha})$, where $\Delta(\sqrt{2\alpha})$ is defined in \cite[(40)]{KLMR}. The equation \eqref{limit-w} follows from the last equality in the proof of \cite[Theorem 2.1 (iii) and (iv)]{KLMR}.

Let ${\cal C}_c(\R)({\cal C}_c^+(\R))$ be the class of all the (nonnegative) continuous functions with compact support.
Let $\cM_{R}(\R)$ be the space of
all the Radon measures on $\R$ equipped with the vague topology, see \cite[p.12]{Kallenberg}.
Recall that for any random measures $\mu_t,\mu\in \cM_{R}(\R)$, $\mu_t{\to} \mu$ in distribution  if and only if for any $f\in{\cal C}_c(\R)$, $\langle f,\mu_t\rangle{\to}\langle f,\mu\rangle$ in distribution,
see \cite[Lemma 4.11]{Kallenberg}.
It follows from \cite[Corollary 4.5]{Kallenberg1} that,
for random measures $\mu_t,\mu\in \cM_{R}(\R)$,
$\mu_t{\to} \mu$ in distribution  is equivalent to $\langle f,\mu_t\rangle{\to}\langle f,\mu\rangle$ in distribution for any $f\in{\cal C}_c^+(\R)$.

For any $z\in \R$ and function $f$ on $\R$, we define  the shift operator $\theta_zf$ by $\theta_zf(y):=f(y+z)$, and for $\mu\in\cM_{R}(\R)$, we define ${\cal T}_z\mu$ by $\int f(y){\cal T}_z\mu(dy):=\int f(y+z)\mu(dy)$. Sometimes, we also write ${\cal T}_z\mu$ as $\mu+z$. We define the rightmost point  $M(\mu)$ of $\mu\in\cM_{R}(\R)$ by
$M(\mu):=\sup\{x: \mu(x,\infty)>0.\}$.
Here we use the convention that $\sup\emptyset=-\infty.$
The supremum $M_t$ of the support of our super-Brownian motion $X_t$ is simply $M(X_t)$.

For any  $\phi\in \mathcal{B}_b^+(\mathbb{R})$, we define
\begin{align}
U_{\phi}(t,x)&:=-\log \E \Big[\exp\Big\{-\int_{\R} \phi(y-x)X_t(dy)\Big\}\Big];\label{def-U}\\
V_{\phi}(t,x)&:=-\log \E\Big[\exp\Big\{-\int_\R \phi(y-x)X_t(dy)\Big\} , M_t\le x\Big];\label{Vphi}\\
V(t,x)&:=-\log \P(M_t\le x).\label{V}
\end{align}
By the spatial homogeneity of $X$,
we have $U_{\phi}(t,x)=u_{\phi}(t,-x)$,
and thus $U_{\phi}(t,x)$ is the unique positive solution to  \eqref{KPPpsi}  with initial condition $U_{\phi}(0,x)=\phi(-x)$.

By the Markov property of $X$, we have
\begin{align*}
 &V_{\phi}(t+r,x)=-\log \E_{\delta_{-x}}\Big[e^{-\int \phi(y)X_{t+r}(dy)} , M_{t+r}\le 0\Big]\\
 &=\lim_{\theta\to\infty}-\log \E_{\delta_{-x}}\Big[ e^{-\langle \phi+\theta\textbf{1}_{(0,\infty)},X_{t+r}\rangle}\Big]
 =\lim_{\theta\to\infty}-\log \E_{\delta_{-x}}\Big[ e^{-\langle u_{\phi+\theta\textbf{1}_{(0,\infty)}}(r),X_t\rangle}\Big]\\
 &=-\log \E_{\delta_{-x}}\Big[ e^{-\int_{\R} V_{\phi}(r,-y),X_t(dy)}\Big]=U_{V_{\phi}(r,-y)}(t,x).
\end{align*}
Thus, for any $r>0$, $(t,x)\to V_{\phi}(t+r,x)$ is a solution to  \eqref{KPPpsi}  with initial condition $V_{\phi}(r,x)$. Thus $V_{\phi}(t,x)$ is a solution to
\eqref{KPPpsi}  with initial condition
\begin{equation}\label{heavy-boundary}
V_{\phi}(0, x)=\left\{
\begin{array}{ll}
  \phi(-x), & x\ge 0;\\
  \infty, & x<0.
   \end{array}
   \right.
\end{equation}

The constants introduced in the next result will be used in the statements of our main results.
\begin{proposition}\label{prop:lim-U}
Assume that $\phi\in \mathcal{B}^+(\R)$ satisfies the following integrability condition at $-\infty$:
\begin{equation}\label{initial-cond1'}
  \int_0^\infty y e^{\sqrt{2\alpha}y}\phi(-y)\,dy<\infty.
\end{equation}
\begin{itemize}
 \item [(1)]
 If $\bf{(H1)}$ and $\bf{(H2)}$ hold, and $\phi$ is bounded, then the limit
$$
C(\phi)=\lim_{r\to\infty}\sqrt{\frac{2}{\pi}}\int_0^\infty U_{\phi}(r,\sqrt{2\alpha}r+y)ye^{\sqrt{2\alpha}y}\,dy
$$
exists,  $C(\theta_z\phi)=C(\phi)e^{\sqrt{2\alpha}z}$ for all $z\in \R$ and
\begin{equation}\label{u(t,sqrt(2)t)}
  \lim_{t\to\infty}\frac{t^{3/2}}{\frac{3}{2\sqrt{2\alpha}}\log t}U_\phi(t,\sqrt{2\alpha}t+x)=C(\phi)e^{-\sqrt{2\alpha}x}, \quad x\in \R.
\end{equation}
If $\phi$ is non-trivial, then $C(\phi)\in (0,\infty)$.

 \item [(2)]
 If $\bf{(H1)}$ and $\bf{(H3)}$ hold, and if there exists $x_0<0$ such that $\phi$ is bounded
 on $(-\infty, x_0]$,
then \eqref{u(t,sqrt(2)t)} holds, and
the limit
$$\tl{C}(\phi):=\lim_{r\to\infty}\sqrt{\frac{2}{\pi}}\int_0^\infty V_{\phi}(r,\sqrt{2\alpha}r+y)ye^{\sqrt{2\alpha}y}\,dy\in(0,\infty)$$
exists and
\begin{equation}\label{v(t,sqrt(2)t)}
 \lim_{t\to\infty}\frac{t^{3/2}}{\frac{3}{2\sqrt{2\alpha}}\log t}V_\phi(t,\sqrt{2\alpha}t+x)=\tl{C}(\phi)e^{-\sqrt{2\alpha}x}, \quad x\in \R.
\end{equation}
\end{itemize}
\end{proposition}

It has been shown in \cite{KLMR} that $\frac{M_t}{t}\to\sqrt{2\alpha}$,
$\P(\cdot|\cS)$-a.s.
Next, we give some large deviation results for $M_t$.

\begin{theorem}\label{Them:large-der}
Under  $\bf{(H1)}$ and $\bf{(H3)}$, the following hold:
\begin{itemize}
  \item[(1)]
\begin{equation}\label{lar-der1}
  \lim_{t\to\infty}\frac{t^{3/2}}{\frac{3}{2\sqrt{2\alpha}}\log t}\P(M_t>\sqrt{2\alpha}t)=\tl{C}_0,
\end{equation}
where $\tl{C}_0$ is the constant $\tl{C}(\phi)$ with $\phi=0$.
  \item[(2)] For any $\delta>0$, the limit
$$\hat{C}(\delta):=\lim_{r\to\infty}\sqrt{\frac{2}{\pi}}\delta e^{-\frac{1}{2}\delta^2 r}\int_0^\infty V(r,\sqrt{2\alpha}r+y)ye^{(\sqrt{2\alpha}+\delta)y}\,dy\in(0,\infty)$$
exists and
\begin{equation}\label{lar-der2}
  \lim_{t\to\infty}\sqrt{t}e^{(\delta^2/2+\sqrt{2\alpha}\delta)t}\P(M_t>(\sqrt{2\alpha}+\delta)t)=\hat{C}(\delta).
\end{equation}
\end{itemize}
\end{theorem}

The analogue of the above results for branching Brownian motions were given in \cite{Chauvin88,Chauvin}.

 \bigskip

In the remainder of this paper, we define
\begin{align}\label{e:m}
m(t):=\sqrt{2\alpha}t-\frac{3}{2\sqrt{2\alpha}}\log t.
\end{align}

\begin{theorem}\label{thrm:tra-wave}
Suppose that $\phi\in\mathcal{B}^+(\R)$ satisfies
the integrability condition \eqref{initial-cond1'} at $-\infty$.
Let $x(\cdot)$ be a function on $\R$ satisfying $\lim_{t\to\infty}x(t)=x\in\R$.
\begin{itemize}
  \item[(1)]
  If $\bf{(H1)}$ and $\bf{(H2)}$ hold, and $\phi$ is bounded, then
  \begin{equation}\label{tra-wave-U}
    \lim_{t\to\infty}U_\phi(t,m(t)+x(t))=-\log \E\left[\exp\{-C(\phi)Z_\infty e^{-\sqrt{2\alpha}x}\}\right].
  \end{equation}
  \item[(2)]
If  $\bf{(H1)}$ and $\bf{(H3)}$ hold, and if there exists $x_0<0$ such that $\phi$ is bounded
 on $(-\infty, x_0]$,  then
  \eqref{tra-wave-U}  holds, and
  $$\lim_{t\to\infty}V_\phi(t,m(t)+x(t))=-\log \E\left[\exp\{-\tilde{C}(\phi)Z_\infty e^{-\sqrt{2\alpha}x}\}\right].$$
\end{itemize}
\end{theorem}

\begin{remark}
In the case when the nonlinear function $f$ satisfies \eqref{f-condition1} and \eqref{f-condition2},
Bramson \cite{Bramson} studied the uniform convergence of solutions of the KPP equation \eqref{KPPequ} to traveling wave solutions.
More precisely, under the integrability condition \eqref{initial-cond1'}  at $-\infty$
and another growth condition of $\phi$ at $+\infty$, he proved that
$u(t, m(t)+x)$ converge uniformly in $x\in\R$, where $u(t,x)$
is the solution of the KPP equation \eqref{KPPequ} with  initial condition
$u(0, x)=\phi(-x)$.
In this paper, our condition on the nonlinear function $-\psi$ is weaker, and we will not
study uniform convergence of solutions of \eqref{KPPpsi} to traveling wave solutions.
\end{remark}

\begin{remark}\label{rm:1.7}
Applying Theorem \ref{thrm:tra-wave}.(2) to $\phi=0$, we get that
\begin{equation}\label{1.1}
  \lim_{t\to\infty}\P(M_t-m(t)\le x)=\E(e^{-\tl{C}_0e^{-\sqrt{2\alpha}x}Z_\infty}),
  \quad x\in \R.
\end{equation}
Using this, one can check that for any $x\in \R$,
$$\lim_{t\to\infty}\P(M_t-m(t)\le x|\cS)=\lim_{t\to\infty}\frac{\P(M_t-m(t)\le x)-\P(\cS^c)}{1-\P(\cS^c)}=\E(e^{-\tl{C}_0e^{-\sqrt{2\alpha}x}Z_\infty}|\cS).$$
Thus, $M_t-m(t)|_{\cS}$ converges in distribution to a random variable $M^*$.
\end{remark}

Let $\mathcal{H}$ be the class of all the nonnegative bounded functions vanishing on $(-\infty,a)$ for some $a\in\R$.
It is clear that the functions in $\mathcal{H}$ satisfy
the integrability condition  \eqref{initial-cond1'} at $-\infty$.
In Lemma \ref{cont-C} below, we will prove that
for any $\phi\in\mathcal{H}$, $C(\lambda\phi)\to0$, $\tl{C}(\lambda\phi)\to\tl{C}_0$ as $\lambda\to0$.
Recall that, for any $t>0$, $\mathcal{E}_t={\cal T}_{-m(t)}X_t$ is the extremal process of $X_t$.
Then $U_\phi(t, m(t))=-\log\E[\exp\{-\langle\phi, \cE_t\rangle\}]$. Using the above theorem, we get that, for any $\phi\in\mathcal{H}$,
\begin{itemize}
  \item [(1)] under $\bf{(H1)}$ and $\bf{(H2)}$,
$\langle\phi,\cE_t\rangle$ converges in distribution;
  \item [(2)] under $\bf{(H1)}$ and $\bf{(H3)}$, $(\langle\phi,\cE_t\rangle,M_t-m(t))|_{\cS}$ jointly converges in distribution.
\end{itemize}
In Theorems \ref{thrm:tra-wave2} and \ref{them:main3}, we will describe these limits.

In Proposition \ref{prop-lim2}, we will prove that, conditioned on $\{M_t>\sqrt{2\alpha}t+z\}$, $X_t-M_t$ converges in distribution to a limit (independent of $z$) denoted by $\Delta$.
Let $\Delta_i,i\ge1,$ be a sequence of independent, identically distributed random variables with the same law as $\Delta$.
Given $Z_\infty,$ let $\sum_{j=1}^\infty \delta_{e_j}$ be a Poisson random measure with intensity $\tilde{C}_0Z_\infty \sqrt{2\alpha}e^{-\sqrt{2\alpha}x}\,dx$. Assume that $\{\Delta_i,i\ge1\}$ and $\sum_{j=1}^\infty \delta_{e_j}$ are independent.

\begin{theorem}\label{thrm:tra-wave2}
Assume that $\bf{(H1)}$ and $\bf{(H2)}$ hold. Then,
as $t\to \infty$,  $\mathcal{E}_t$ converges in law to a random Radon measure $\mathcal{E}_\infty$ with Laplace transform
\begin{equation}\label{laplace-E-M}
  \E\Big[\exp\Big\{-\int \phi(y)\mathcal{E}_\infty(dy)\Big\}\Big]
  =\E\Big[\exp\left\{-C(\phi)Z_\infty\right\}\Big], \quad \phi\in\mathcal H.
\end{equation}
Moreover, if, in addition, $\bf{(H3)}$ holds, then
$$\mathcal{E}_\infty\overset{d}{=}\sum_{j}{\cal T}_{e_j}\Delta_j.$$
\end{theorem}

For any $t>0$, we define $\cE^*_t=X_t-m(t)-\frac{1}{\sqrt{2}}\log Z_\infty$. Then we have the following result.

\begin{theorem}\label{them:main3}
Assume that $\bf{(H1)}$ and $\bf{(H3)}$ hold.
Conditioned on  $\cS$,
$(\cE^*_t,Z_t)$ converges jointly in distribution to $(\cE^*_\infty,Z_\infty^*)$, where $Z_\infty^*$ has the same law as $Z_\infty$ conditioned on $\cS$, $\cE^*_\infty$ and $Z_\infty^*$ are independent,
and the Laplace transform of $\cE^*_\infty$ is given by
\begin{equation}\label{laplace1}
  \E[\exp\{-\langle \phi,\cE^*_\infty\rangle\}]=\exp\{-C(\phi)\},
  \quad \phi\in\mathcal{C}_c^+(\R).
\end{equation}
Moreover,
\begin{equation}\label{decom:E*}
  \cE^*_\infty\overset{d}{=}\sum_{j}{\cal T}_{e_j'}\Delta_j,
\end{equation}
where $\sum_{j}\delta_{e_j'}$ is a Poisson random measure with intensity measure $\tl{C}_0\sqrt{2\alpha}e^{-\sqrt{2\alpha}x}\,dx$, which is independent of $\{\Delta_j,j\ge 1\}$.
\end{theorem}

Following the arguments in Hu and Shi \cite{HS} for branching random walks (see also  Roberts \cite{Robert} for branching Brownian motions), together with
Lemma \ref{lem-h} and \ref{lem:com2} below,
one can prove some
almost sure fluctuation results $M_t$. We will not pursue this in this paper.

Note that it suffices to  prove the above results for the case $\alpha=1$ and $\lambda^*=1$.
For the  general case, let
$v(t,x)=\frac{1}{\lambda^*}u(\alpha^{-1}t,\alpha^{-1/2}x)$. If $u$ is a non-negtaive solution of
\eqref{KPPpsi}, then, $v$ is a non-negative solution of \eqref{KPPpsi} with $\psi$ replaced by
where  $\psi^*(x)=\frac{\psi(\lambda^* x)}{\alpha\lambda^*}$.
It is clear that $-\psi^*$ satisfies
condition \eqref{f-condition1}.
Therefore, in the remainder of this paper,  we assume that $\alpha=1$ and $\lambda^*=1$.

The rest of the paper is organized as follows. In Section 2, we generalize some results in \cite{Bramson}
to the case when the nonlinear term satisfies a weaker condition and to general initial conditions.
In Section 3.1, we give the proofs of the large deviation results,
 including Proposition 1.1 and Theorem 1.2.
In Section 3.2, we study the convergence of the extremal process.
 In Section 4, we give the proof of Lemma \ref{initial-cond-V}.

\section{Some results on the KPP equation \eqref{KPPpsi}}
It follows from the Feyman-Kac formula that, if $u$ is a non-negative solution to
\eqref{KPPpsi}, then, for any $0\le r<t$,
\begin{equation}\label{F-K-u}
  u(t,x)=\mE_x\left[u(r,B_{t-r})\exp\Big\{\int_0^{t-r} k(u(t-s,B_s))\,ds\Big\}\right],
\end{equation}
where $k(\lambda)=\frac{-\psi(\lambda)}{\lambda}$. Recall that we always assume that $\alpha=1$ and $\lambda^*=1$.
Note that $k(\lambda)$ is decreasing and $k(\lambda)\le 1$ for all $\lambda>0$.
We first give some basic results on non-negative solutions $u$ of the KPP equation \eqref{KPPpsi}
with initial conditions $u(0, \cdot)$ not necessarily bounded between 0 and 1.

\begin{lemma}\label{lem:initial-cond}
Assume that $u(t,x)$ is a solution to \eqref{KPPpsi} with initial condition $u(0,\cdot)\in\mathcal{B}_b^+(\R)$ satisfying
the following integrability condition at $\infty$:
\begin{equation}\label{initial-cond1}
  \int_0^\infty y e^{\sqrt{2}y}u(0,y)\,dy<\infty.
\end{equation}
Then for any $t>0$,
$u(t,\cdot)$ also satisfies  \eqref{initial-cond1} and
\begin{equation}\label{3.4}
  \int_0^\infty u(r,\sqrt{2}r+y)ye^{\sqrt{2}y}\,dy<\infty.
\end{equation}
\end{lemma}
\textbf{Proof:}
By \eqref{F-K-u} with $r=0$, we have that $u(t,x)\le e^t\mE_x(u(0,B_t))$. Thus, it suffices to show that $\mE_x(u(0,B_t))$ satisfies  \eqref{initial-cond1}. Note that
\begin{align*}
  \int_0^\infty \mE_y(u(0,B_t)) ye^{\sqrt{2}y}\,dy & =\int_{\R}\int_0^\infty u(0,x+y) ye^{\sqrt{2}y}\,dy \mP_0(B_t\in dx) \\
   & =\int_{\R}\int_x^\infty u(0,y) (y-x)e^{\sqrt{2}(y-x)}\,dy \mP_0(B_t\in dx).
\end{align*}
If $x>0$, we have
$$\int_x^\infty u(0,y) (y-x)e^{\sqrt{2}(y-x)}\,dy\le \int_0^\infty u(0,y) ye^{\sqrt{2}y}\,dy,$$
and if $x\le 0$, we have
\begin{align*}
  &\int_x^\infty u(0,y) (y-x)e^{\sqrt{2}(y-x)}\,dy\\
  &  =\int_x^0 u(0,y) (y-x)e^{\sqrt{2}(y-x)}\,dy+\int_0^\infty u(0,y) (y-x)e^{\sqrt{2}(y-x)}\,dy \\
  & \le \|u_0\|_\infty |x|e^{-x\sqrt{2}}+\int_0^\infty u(0,y) ye^{\sqrt{2}y}\,dye^{-x\sqrt{2}}+|x|e^{-x\sqrt{2}}\int_0^\infty u(0,y)e^{\sqrt{2}y}\,dy\\
  &\le c(1+|x|)e^{-x\sqrt{2}}.
\end{align*}
Thus,
$$\int_0^\infty y e^{\sqrt{2}y}u(t,y)\,dy\le e^t\int_0^\infty \mE_y(u(0,B_t)) ye^{\sqrt{2}y}\,dy <\infty.$$
Using a similar argument, we can get \eqref{3.4}.
The proof is now complete.
\hfill$\Box$

\begin{lemma}{\bf (Maximum principle)}\label{l:mp}
Let $v_1(t,x)$ and $v_2(t,x)$ be non-negative functions satisfying
$$\frac{\partial}{\partial t}v_{2}(t,x)-\frac{1}{2}\frac{\partial^2}{\partial x^2}v_{2}(t,x)+\psi(v_2(t,x))\ge \frac{\partial}{\partial t}v_{1}(t,x)-\frac{1}{2}\frac{\partial^2}{\partial x^2}v_{1}(t,x)+\psi(v_1(t,x)), \quad t>0, x\in(a,b),$$
and
$$v_1(0,x)\le v_2(0,x),\quad x\in(a,b),$$
where $-\infty\le a<b\le \infty$. Moreover, if $a>-\infty$, we assume $v_1(t,a)\le v_2(t,a)$ for all $t>0$, and if $b<\infty$, we assume $v_1(t,b)\le v_2(t,b)$ for all $t>0$.
Then we have that
$$v_1(t,x)\le v_2(t,x),\quad t>0, x\in(a,b).$$
\end{lemma}
\textbf{Proof:}
The proof is a slight modification of the proof of \cite[Proposition 3.1]{Bramson}, using \cite[Theorem 3.4]{PW}. See also the proof of \cite[Proposition 6.4]{Bovier}. We omit the details.
\hfill$\Box$

\begin{lemma}\label{lem:compare}
Assume that $u_1$, $u_2$ and $u_3$
are solutions to  \eqref{KPPpsi} with non-negative bounded initial conditions.
\begin{enumerate}
  \item[(1)] If for some $c>1$, $u_1(0, \cdot)\le cu_2(0, \cdot)$, then $u_1(t,x)\le cu_2(t,x)$
  for all $(t, x)\in (0, \infty)\times \R$.
  \item[(2)] If $u_3(0, \cdot)\le u_1(0, \cdot)+u_2(0, \cdot)$, then $u_3(t,x)\le u_1(t,x)+u_2(t,x)$ for all $(t, x)\in (0, \infty)\times \R$.
\end{enumerate}
\end{lemma}
\textbf{Proof:} (1) Let $v(t,x)=cu_2(t,x)$. Then
$$v_t-\frac{1}{2}v_{xx}=-c\psi(u_2)\ge -\psi(v),$$
here we used the fact that $\psi'(\lambda)$ is increasing.
Applying the maximum principle in Lemma \ref{l:mp}, we get that
$$
u_1(t,x)\le cu_2(t,x), \quad (t, x)\in (0, \infty)\times \R.
$$

(2) Let $v(t,x):=u_1(t,x)+u_2(t,x).$
  Since $\psi'(\lambda)$ is increasing, for any $\theta>0$, the function $\lambda\to\psi(\lambda+\theta)-\psi(\lambda)-\psi(\theta) $ is increasing, which implies that $\psi(\lambda+\theta)\ge \psi(\lambda)+\psi(\theta)$. Then
  $$v_t-\frac{1}{2}v_{xx}=-\psi(u_1)-\psi(u_2)\ge -\psi(v).$$
Applying the maximum principle in Lemma \ref{l:mp}, we get that
$u_3(t,x)\le u_1(t,x)+u_2(t,x)$ for all $(t, x)\in (0, \infty)\times \R$.
\hfill$\Box$

For any $\lambda< 1$ and $y>e^{2+\gamma}$,
one can easily check that $(\lambda y\wedge1)\le  |\log \lambda|^{-2-\gamma}(\log y)^{2+\gamma}.$
Thus, for any $\lambda< 1$,
\begin{align*}
  &0\le 1+\psi'(\lambda)=2\beta\lambda+\int_0^\infty y(1-e^{-\lambda y})n(dy)\nonumber\\
  &\le \Big(2\beta+\int_0^{e^{2+\gamma}} y^2 n(dy)\Big)\lambda+|\log \lambda|^{-(2+\gamma)}\int_{e^{2+\gamma}}^\infty y(\log y)^{2+\gamma}n(dy)\le c_1 |\log{\lambda}|^{-(2+\gamma)},
\end{align*}
where $\gamma$ is the constant in {\bf (H1)} and $c_1>0$ is a constant. Thus {\bf (H1)} implies
\begin{equation}\label{psi-condition2}
 1+\psi'(\lambda)\le c_1|\log{\lambda}|^{-(2+\gamma)}\quad\mbox{ for }\lambda<1.
\end{equation}
Since $\psi'(\lambda)$ is increasing, we have $-k(\lambda)=\psi(\lambda)/\lambda\le \psi'(\lambda)$. Thus
\begin{equation}\label{est-k}
    0\le 1-k(\lambda) \le c_1 |\log{\lambda}|^{-(2+\gamma)},\quad\mbox{ for }\lambda<1.
\end{equation}

In the remainder of this section, we will generalize \cite[Proposition 8.3]{Bramson} to
non-negative solutions of \eqref{KPPpsi} with initial conditions not necessarily bounded between
0 and 1.
The main idea of the proof  is similar to that of \cite{Bramson}.
For the KPP equation \eqref{KPPpsi}, $-\psi$ plays the role of $f$ in \cite{Bramson}. Condition \eqref{f-condition2} is translated to  the following condition on $\psi$:
\begin{equation}\label{psi-condition}
1+\psi'(\lambda)=O(\lambda^\rho)\quad\mbox{as }\lambda\to0, \mbox{ for some }\rho>0.
\end{equation}
However, many results in \cite{Bramson} still hold under the weaker condition \eqref{psi-condition2}.  We will clearly spell out the reason when we apply results from \cite{Bramson}
under this weaker condition.

In the rest of this section, we use $u(t,x)$ to denote the solution to
\eqref{KPPpsi} with initial condition $u(0,\cdot)\in \mathcal{B}_b^+(\R)$.
Let $\tilde{u}(t,x)$ be the solution to \eqref{KPPpsi} with
$\tilde{u}(0,\cdot)=u(0,\cdot)\wedge 1$.
Then, it is clear that
$$\tilde{u}(0,x)\le u(0,x) \le s_u\tilde{u}(0,x), \quad x\in \R,
$$
where $s_u=\sup_{x}u(0,x)\vee 1.$
It follows from Lemma \ref{lem:compare} that
\begin{equation}\label{com-u}
  \tilde{u}(t,x)\le u(t,x)\le s_u\tilde{u}(t,x), \quad (t,  x)\in (0, \infty)\times \R.
\end{equation}
Since $\tilde{u}(t,x)\in[0,1]$, we have
\begin{equation}\label{bound-u3}
  0\le u(t,x)\le s_u.
\end{equation}
\eqref{com-u} and \eqref{bound-u3} will play key roles later in this paper, they allow us to ``scale'' solutions to
KPP equations for Bramson's results to carry over.
Let $\tilde{m}(t)$ be the median of $\tilde{u}$, that is
$$\tl{m}(t):=\sup\{x:\tl{u}(t,x)\ge 1/2\}.$$
It was proved in \cite[(3.22')]{Bramson}, without using condition \eqref{f-condition2}
(equivalently, \eqref{psi-condition}), that
$$\tl{m}(t)/t\to\sqrt{2},\quad t\to\infty.$$

Now we recall some notation from \cite{Bramson}, see \cite[(6.11)--(6.14), (7.6)--(7.9), (7.42), (7.44)]{Bramson}.
In the list of notation below, $\delta\in(\frac{1}{2+\gamma},1/2)$,  $r>1$ and $t>3r$.
\begin{itemize}
\item If $L$ is a function on $[0,t]$, define
  $$\theta_{r,t}\circ L(s):=\left\{
                              \begin{array}{ll}
                                L(s+s^\delta)+4s^\delta, & \hbox{$r\le s\le t/2$;} \\
                                L(s+(t-s)^\delta)+4(t-s)^\delta, & \hbox{$t/2\le s\le t-2r$;} \\
                                L(s), & \hbox{otherwise.}
                              \end{array}
                            \right.
  $$
  \item The inverse of $\theta_{r,t}$ is defined by
  $$
  \theta_{r,t}^{-1}\circ L:=\inf\{l:\theta_{r,t}\circ l\ge L\},
  $$
  that is
  \begin{equation}\label{def:theta-1}
    \theta_{r,t}^{-1}\circ L(s)=\left\{
      \begin{array}{ll}
       -\infty, & \hbox{$r\le s<r+r^\delta$;}\\
        L(u)-4u^\delta, & r+r^\delta\le s\le t/2+(t/2)^\delta; \\
        L(u)-4(t-u)^\delta, & t/2+(t/2)^\delta\le s\le t-2r;\\
        (L(u)-4(t-u)^\delta)\vee L(s),&t-2r<s<t-2r+(2r)^{\delta};\\
        L(s), &\hbox{otherwise,}
      \end{array}
    \right.
  \end{equation}
  where for $r+r^\delta\le s\le t/2+(t/2)^\delta$, $u$ is determined by $s=u+u^\delta$;
  for $t/2+(t/2)^\delta\le s\le t-2r+(2r)^{\delta}$, $u$ is determined by $s=u+(t-u)^\delta$.
  \item $L_{r,t}(s):={\tl m}(s)-\frac{s}{t}{\tl m}(t)+\frac{t-s}{t}\log r,\qquad 0\le s\le t.$
  \item $\underline{\cL}_{r,t}(s):=\theta_{r,t}^{-1}\circ L_{r,t}(s).$
  \item $\overline{\cL}_{r,t}(s):=\theta_{r,t}\circ L_{r,t}(s)\vee \underline{\cL}_{r,t}(s)\vee L_{r,t}(s).$
  \item For any $x$, define
  \begin{equation}\label{def-oM}
    \ocM^x(s):=\left\{
    \begin{array}{ll}
      \ocL(s)+\frac{s}{t}{\tl m}(t)-\frac{t-s}{t}\log r, & \hbox{$0\le s\le t-2r$;} \\
      \frac{x+{\tl m}(t)}{2}, & \hbox{$t-2r<s\le t$.}
    \end{array}
  \right.
\end{equation}
  \item
  \begin{equation}\label{def-uM}
    \ucM'(s):=\left\{
                     \begin{array}{ll}
                       \ucL(s)+\frac{s}{t}{\tl m}(t)-\frac{t-s}{t}\log r, & \hbox{$r+r^\delta\le s\le t-2r$;} \\
                       -\infty, & \hbox{otherwise.}
                     \end{array}
                   \right.
\end{equation}
\item
\begin{equation}\label{e:nrt}
n_{r,t}(s)=\sqrt{2}r+\frac{s-r}{t-r}(m(t)-\sqrt{2}r),\quad s\in[r,t].
\end{equation}
\end{itemize}

The following lemma says that \cite[Proposition 7.2]{Bramson} still holds  without condition \eqref{psi-condition}.
\begin{lemma}\label{prop7.2}
Assume that (H1) holds. Let $u(t,x)$ be a  solution to \eqref{KPPpsi} with initial condition $u(0,\cdot)\in\mathcal{B}_b^+(\R)$ satisfying \eqref{initial-cond1}.
For  all $t>3r>0$, and continuous function $x(s)$ with $x(s)>\ocM^x(t-s)$ in $[2r,t-r]$, we have that
\begin{align}\label{bound-int-k}
  e^{3r-t}\int_{2r}^{t-r}k(u(t-s,x(s)))\,ds\to1,\quad r\to\infty,
\end{align}
uniformly in $t$.
\end{lemma}
{\bf Proof:}
First note that the proofs  of \cite[(7.16) and (7.18)]{Bramson} did not use
\eqref{f-condition2} (equivalently \eqref{psi-condition}).
Thus
there exists a constant $C>0$ such that for $r$ large enough, $s\in[r,t-2r]$ and $y>{\tl m}(s+(s\wedge(t-s))^\delta)$,
$$
  \tilde u(s,y)\le Ce^{-(s\wedge(t-s))^\delta}.
$$
It follows immediately from the key inequality \eqref{com-u} that
\begin{equation}\label{bound-u1}
  u(s,y)\le c_2e^{-(s\wedge(t-s))^\delta}.
\end{equation}
For $r$ large enough and $s\in[r,t-2r]$, by the definition of $\ocM^x$, we have
\begin{align*}
  &\ocM^x(s)\ge  \theta_{r,t}\circ L_{r,t}(s)+\frac{s}{t}{\tl m}(t)-\frac{t-s}{t}\log r\\
  &=L_{r,t}(s+(s\wedge(t-s))^\delta)+4(s\wedge(t-s))^\delta+\frac{s}{t}{\tl m}(t)-\frac{t-s}{t}\log r\\
  &=\tl m(s+(s\wedge(t-s))^\delta)+(s\wedge(t-s))^\delta(4-\log r/t-\tl{m}(t)/t)\ge \tl m(s+(s\wedge(t-s))^\delta),
\end{align*}
where in the last inequality, we use the fact that $4-\log r/t-\tl{m}(t)/t\ge 4-\log t/t-\tl{m}(t)/t\to 4-\sqrt{2}>0$ as $t\to\infty$.
Thus, by \eqref{bound-u1}, for $r$ large enough, we have
\begin{align}\label{3.2'}
 &e^{t-3r}\ge\exp\{\int_{2r}^{t-r}k(u(t-s,x(s)))\,ds\} \nonumber\\
  \ge &e^{t-3r}\exp\Big\{-\int_{2r}^{t-r}(1-k(c_2e^{-(s\wedge(t-s))^\delta}))\,ds\Big\} \nonumber\\
  \ge &e^{t-3r}\exp\Big\{-2\int_{r}^{\infty}(1-k(c_2e^{-s^\delta}))\,ds\Big\}.
\end{align}
For $\delta>\frac{1}{2+\gamma}$, by \eqref{est-k},
\begin{align*}
  \int_{r}^{\infty}(1-k(c_2e^{-s^\delta}))\,ds & \le \int_{r}^{\infty}c_1|s^{\delta}-\log c_2|^{-(2+\gamma)}\,ds\to 0, \quad r\to\infty.
\end{align*}
Thus, the desired result follows immediately.
\hfill$\Box$

The lemma above implies that, under {\bf (H1)}, (7.12) in  \cite[Proposition 7.2]{Bramson} is valid for $\tl{u}(t,x)$. Since in the proofs of \cite[Proposition 8.1,  Corollary 1 on p. 125, Proposition 8.2, Corollary 1 on p. 130 and Corollary 2 on p. 133]{Bramson},
only \cite[(7.12)]{Bramson} was used, these results hold for $\tl{u}(t,x)$
under {\bf (H1)}.
Thus,
\begin{equation}\label{equa-mt}
  \tl{m}(t)=m(t)+O(1),
\end{equation}
where
\begin{equation}\label{def-mt}
  m(t)=\sqrt{2}t-\frac{3}{2\sqrt{2}}\log t.
\end{equation}

\begin{proposition}\label{prop:u-psi}
Assume that $u(t,x)$ is a solution to \eqref{KPPpsi} with initial condition $u(0,\cdot)\in\mathcal{B}_b^+(\R)$ satisfying \eqref{initial-cond1}.
Define, for $0\le r\le t$,
\begin{equation}\label{e:Psi}
\Psi(r,t,x):=\frac{e^{-\sqrt{2}(x-\sqrt{2}t)}}{\sqrt{2\pi(t-r)}}\int_0^\infty u(r,\sqrt{2}r+y)e^{\sqrt{2}y}e^{-\frac{(x-\sqrt{2}t-y)^2}{2(t-r)}}(1-e^{-2(x-m(t))y/(t-r)})\,dy,
\end{equation}
where $m(t)$ is defined in \eqref{def-mt}.
Then for $r$ large enough, $t\ge8r$ and $x\ge m(t)+9r$,
\begin{equation}\label{comp-u-psi}
  \gamma(r)^{-1}\Psi(r,t,x)\le u(t,x)\le \gamma(r)\Psi(r,t,x),
\end{equation}
where $\gamma(r)\downarrow 1$ as $r\to\infty$.
\end{proposition}

To prove the proposition above, we need the following lemma whose proof is similar to that of \cite[(8.62)]{Bramson}.
Let $(B_{x,y}^t, P)$ be a Brownian bridge starting from $x$ and ending at $y$ at time $t$, and $E$ be the expectation with respect to $P$.

\begin{lemma}\label{lem:bound-u}
Assume that $u(t,x)$ solves the KPP equation \eqref{KPPpsi} with initial condition $u(0,\cdot)\in\mathcal{B}_b^+(\R)$ satisfying \eqref{initial-cond1}.
Then for large $r$, $t>8r$ and $x\ge \tilde{m}(t)+8r$,
\begin{align}\label{low-bound-u}
  &u(t,x)\ge \psi_1(r,t,x)\nonumber\\
  &:=C_1(r)e^{t-r}\int_{-\infty}^\infty u(r,y)\frac{e^{-\frac{(x-y)^2}{2(t-r)}}}{\sqrt{2\pi (t-r)}}P\Big[B_{x,y}^{t-r}(s)>\bar{\mathcal{M}}_{r,t}^x(t-s),s\in[0,t-r]\Big]\,dy
\end{align}
and
\begin{align}\label{up-bound-u}
  &u(t,x)\le \psi_2(r,t,x)\nonumber\\
  &:=C_2(r)e^{t-r}\int_{-\infty}^\infty u(r,y)\frac{e^{-\frac{(x-y)^2}{2(t-r)}}}{\sqrt{2\pi (t-r)}}P\Big[B_{x,y}^{t-r}(s)>\underline{\mathcal{M}}_{r,t}'(t-s),s\in[0,t-r]\Big]\,dy
\end{align}
with $C_1(r)\to 1,$ $C_2(r)\to 1,$ as $r\to\infty.$
Moreover,
\begin{equation}\label{com-psi}
  1\le \frac{\psi_2(r,t,x)}{\psi_1(r,t,x)}\le \gamma(r),
\end{equation}
with $\gamma(r)\downarrow 1,$ as $r\to\infty$.
\end{lemma}
\textbf{Proof:}
Let $${\cal A}:=\Big\{B_{x,y}^{t-r}(s)>\overline{\mathcal{M}}_{r,t}^x(t-s),s\in[0,t-r]\Big\}.$$
It follows from \eqref{F-K-u} that
\begin{align}\label{formula-u2}
  u(t,x)&=\int_{-\infty}^\infty u(r,y)\frac{e^{-\frac{(x-y)^2}{2(t-r)}}}{\sqrt{2\pi (t-r)}}E \Big[\exp\Big\{\int_{0}^{t-r}k(u(t-s,B_{x,y}^{t-r}(s)))\,ds\Big\}\Big]\,dy\nonumber\\
  &\ge \int_{-\infty}^\infty u(r,y)\frac{e^{-\frac{(x-y)^2}{2(t-r)}}}{\sqrt{2\pi (t-r)}}E \Big[\exp\Big\{\int_{0}^{t-r}k(u(t-s,B_{x,y}^{t-r}(s)))\,ds\Big\},{\cal A}\Big]\,dy.
\end{align}
For $r$ large enough, $t>8r$, $s\in[0,2r]$ and $x\ge \tilde{m}(t)+8r$, it holds that
$$
\ocM^x(t-s)=(x+\tl{m}(t))/2\ge \tl{m}(t)+4r\ge \tl{m}(t-s)+r,
$$
where in the last inequality, we used the fact that
$\tl{m}(t)-\tl{m}(t-s)=m(t)-m(t-s)+O(1)$ is bounded from below, because $m(t)$ is increasing on $t\ge1$.
Thus, by \eqref{com-u} first, and then
applying \cite[Proposition 8.2]{Bramson} to $\tl{u}$,
we get that on ${\cal A}$,
$$u(t-s,B_{x,y}^{t-r}(s))\le s_u \tl{u}(t-s,B_{x,y}^{t-r}(s))\le c_3re^{-\sqrt{2}r}.$$
It follows from \eqref{est-k} that for $r$ large enough,
\begin{align}\label{3.1}
   &E \Big[\exp\{\int_0^{2r}k(u(t-s,B_{x,y}^{t-r}(s)))\,ds\},{\cal A}\Big] \nonumber  \\
& \ge e^{2r}\exp\Big\{-2rc_1|\log (c_3re^{-\sqrt{2}r})|^{-(2+\gamma)})\Big\}P ({\cal A}).
\end{align}
Note that $2r|\log (c_3re^{-\sqrt{2}r})|^{-(2+\gamma)}\to 0$ as  $r\to\infty$.

By Lemma \ref{prop7.2}, we have
\begin{align}\label{3.2}
 &E \Big[\exp\{\int_{2r}^{t-r}k(u(t-s,B_{x,y}^{t-r}(s)))\,ds\},{\cal A}\Big] \nonumber\\
  \ge &e^{t-3r}\int_{r}^{\infty}c_1|s^{\delta}-\log c_2|^{-(2+\gamma)}\,ds P({\cal A}).
\end{align}
Combining \eqref{formula-u2}--\eqref{3.2}, we immediately get \eqref{low-bound-u}.
The proof of \eqref{up-bound-u} is similar to that of \cite[Proposition 8.3 (b)]{Bramson} and the proof of \eqref{com-psi} is similar to that of \cite[(8.62)]{Bramson}.
Here we omit the details.
\hfill$\Box$

\noindent\textbf{Proof of Proposition \ref{prop:u-psi}:}
Recall that $n_{r,t}(\cdot)$ is defined by \eqref{e:nrt}.
First, we claim that, for $s\in[r,t]$,
\begin{equation}\label{com-M-n}
  \ucM'(s)\le n_{r,t}(s)\le \ocM^x(s).
\end{equation}
It has been proved in \cite[Lemma 2.2]{Bramson} that for $y>\sqrt{2}r$ and $x>m(t)$,
\begin{align*}
 &P
  \Big[B_{x,y}^{t-r}(s)>n_{r,t}(t-s),s\in[0,t-r]\Big]\\
 =&P
  \Big[B_{0,0}^{t-r}(s)>-\frac{s}{t-r}(y-\sqrt{2}r)-\frac{t-r-s}{t-r}(x-m(t)),s\in[0,t-r]\Big]\\
  =&1-\exp\Big\{-\frac{2(x-m(t))(y-\sqrt{2}r)}{t-r}\Big\},
\end{align*}
and for $y\le \sqrt{2}r$, $P\Big[B_{x,y}^{t-r}(s)>n_{r,t}(t-s),s\in[0,t-r]\Big]=0$.
Thus, combining  Lemma \ref{lem:bound-u} and \eqref{com-M-n}, the desired result follows immediately.

Now we prove the claim. For $r$ large enough, $s\in[r+r^\delta,t/2]$ and $u$ determined by $s=u+u^\delta$,
\begin{align*}
   \ucM'(s)&=L_{r,t}(u)-4u^\delta+\frac{s}{t}\tl{m}(t)-\frac{t-s}{t}\log r \\
   & =\tl{m}(u)-u^{\delta}(4-\frac{\log r+\tl{m}(t)}{t})\le m(u)\le m(s),
\end{align*}
where in the last inequality we used that fact $m(t)$ is increasing for $t$ large enough.
Similarly, for $s\in[t/2,t-2r]$, $\ucM'(s)\le m(s)$.
Thus, for all $s\in[0,t]$, $\ucM's)\le m(s)$.

By the definition of $n_{r,t}(s)$, for $r$ large enough,
$$n_{r,t}(s)=\sqrt{2}r+\frac{s-r}{t-r}(m(t)-\sqrt{2}r)=\sqrt{2}s-\frac{s-r}{t-r}\frac{3}{2\sqrt{2}}\log t \ge m(s),$$
where we used the fact that for large $r$, $t\to \log t/(t-r)$ is decreasing.
Thus, we get that $\ucM'(s)\le n_{r,t}(s)$, $s\in[r,t]$.

Now we deal with $\ocM^x(s)$. For $r$ large enough, $s\in[r,t/2]$,
$$\ocM^x(s)\ge m(s+s^{\delta}))\ge \sqrt{2}s\ge n_{r,t}(s).$$
For $r$ large enough, $s\in[t/2,t-2r]$,
$$\ocM^x(s)\ge m(s+(t-s)^{\delta})\ge \sqrt{2}s+\sqrt{2}(t-s)^{\delta}-\frac{3}{2\sqrt{2}}\log t\ge n_{r,t}(s).$$
For $r$ large enough, $s\in[t-2r,t]$ and $x\ge m(t)+9r$,
$$\ocM^x(s)=\frac{\tl{m}(t)+x}{2}\ge m(t)\ge n_{r,t}(s),$$
where the last inequality follows from the fact that, for $r$ large enough, $\sqrt{2}r\le m(t)$.
The proof is now complete.
\hfill$\Box$

\section{Proof of main results}
We first give a useful lemma. The proof of this lemma will be given in Section 5.
\begin{lemma}\label{initial-cond-V}
Assume that ${\bf (H1)}$ and ${\bf (H3)}$ hold.
Then, for any $t>0$ and $\theta>0$, we have that $V(t,\cdot)\in \mathcal{B}_b^+(\R)$,  and
$$\int_0^\infty V(t,x)xe^{\theta x}\,dx<\infty.$$
\end{lemma}
\smallskip
\begin{corollary}\label{initial-cond-V2}
Assume that ${\bf (H1)}$ and ${\bf (H3)}$ hold.
If $\phi\in\mathcal{B}^+(\R)$ satisfies the integrability condition \eqref{initial-cond1'}
at $-\infty$ and if there exists $x_0<0$ such that $\phi$ is bounded on $(-\infty, x_0]$, then
for any $t>0$, we have that $U_\phi(t,\cdot)$ and $V_\phi(t,\cdot)$ are bounded
functions satisfying \eqref{initial-cond1}.
\end{corollary}
\noindent\textbf{Proof:}
First, we assume that $\phi\in\mathcal{B}_b^+(\R)$ and satisfies the integrability condition
\eqref{initial-cond1'} at $-\infty$.
Note that $U_\phi(t,x)$ satisfies
  the KPP equation \eqref{KPPpsi}
with $U_\phi(0,x)=\phi(-x)\in\mathcal{B}_b^+(\R)$ satisfying the integrability condition
\eqref{initial-cond1} at $\infty$,
thus by Lemma \ref{lem:initial-cond}, $U_\phi(t,x)$ is  bounded and satisfies \eqref{initial-cond1}.
For $V_\phi(t,x)$, it is clear that, by Lemma \ref{lem:compare}.(2),
\begin{align*}
&  V_\phi(t,x) =\lim_{\theta\to\infty} U_{\phi+\theta{\bf 1}_{(0,\infty)}}(t,x)\\
&\le  U_{\phi}(t,x)+\lim_{\theta\to\infty} U_{\theta1_{(0,\infty)}}(t,x) \le  U_{\phi}(t,x)+V(t,x),
\end{align*}
where $V$ is defined in \eqref{V}.
By Lemma \ref{initial-cond-V},  $V(t,x)$ is bounded and satisfies the integrability condition
\eqref{initial-cond1} at $\infty$.
Thus $V_\phi$ is a bounded function satisfying the integrability condition
\eqref{initial-cond1} at $\infty$,
when $\phi\in\mathcal{B}_b^+(\R)$ and satisfies the integrability condition
\eqref{initial-cond1'} at $-\infty$.

Now we assume that $\phi\in\mathcal{B}^+(\R)$ satisfies the integrability condition \eqref{initial-cond1'}
at $-\infty$ and that there exists $x_0<0$ such that $\phi$ is bounded on $(-\infty, x_0]$. Let $\tl{\phi}(x):=\phi(x)1_{x\le x_0}\in\mathcal{B}_b^+(\R)$.
Note that
$$V_{\phi}(t,x)\le -\log \E\Big[\exp\Big\{-\int_\R {\tl \phi}(y-x)X_t(dy)\Big\} ,
 M_t\le x+x_0\Big]
=V_{\tl{\phi}(\cdot+x_0)}(t, x+x_0).
$$
Thus $V_\phi$ is a bounded function satisfying the integrability condition
\eqref{initial-cond1} at $\infty$.  Since $U_\phi(t,x)\le V_\phi(t,x)$, $U_\phi$
is also a bounded function satisfying the integrability condition
\eqref{initial-cond1} at $\infty$.
\hfill $\Box$

\subsection{Large deviation results}

\noindent\textbf{Proof of Proposition \ref{prop:lim-U}:}
(1)
Let $\Psi$ be defined by \eqref{e:Psi} with $u$ replaced by $U_\phi$.
We claim that
\begin{equation}\label{lim-psi1}
  \lim_{t\to\infty}e^{\sqrt{2}x}\frac{t^{3/2}}{\frac{3}{2\sqrt{2}}\log t}\Psi(r,t,\sqrt{2}t+x)=\sqrt{\frac{2}{\pi}}\int_0^\infty U_{\phi}(r,\sqrt{2}r+y)ye^{\sqrt{2}y}\,dy:=C(\phi, r).
\end{equation}
Note that, for $r>1$ and $t\ge 8r$,
\begin{align*}
  &e^{\sqrt{2}x}\frac{t^{3/2}}{\frac{3}{2\sqrt{2}}\log t}\frac{e^{-\sqrt{2}x}}{\sqrt{2\pi(t-r)}}U_\phi(r,\sqrt{2}r+y)e^{\sqrt{2}y}e^{-\frac{(x-y)^2}{2(t-r)}}\left(1-e^{-2(x+\frac{3}{2\sqrt{2}}\log t)y/(t-r)}\right)\\
  &\le c(1+|x|)U_\phi(r,\sqrt{2}r+y)ye^{\sqrt{2}y}.
\end{align*}
By \eqref{3.4}, the right hand side of the inequality above is integrable. So by the dominated convergence theorem, we get that the claim is true and that $C(\phi, r)\in(0,\infty).$

Since $U_\phi$ is the solution to \eqref{KPPpsi} with $U_{\phi}(0,x)=\phi(-x)$ satisfying
the integrability condition  \eqref{initial-cond1} at $\infty$,
by Proposition \ref{prop:u-psi}, for $r$ large enough, $t\ge 8r$ and $x\ge -\frac{3}{2\sqrt{2}}\log t+9r$,
$$
  \gamma(r)^{-1}\Psi(r,t,\sqrt{2}t+x)\le U_\phi(t,\sqrt{2}t+x)\le \gamma(r)\Psi(r,t,\sqrt{2}t+x).
$$
Thus, by \eqref{lim-psi1}, we have
\begin{align}\label{3.6}
  &\gamma(r)^{-1}C(\phi,r)\le \liminf_{t\to\infty} e^{\sqrt{2}x}\frac{t^{3/2}}{\frac{3}{2\sqrt{2}}\log t}U_\phi(t,\sqrt{2}t+x)\nonumber\\
  &\le \limsup_{t\to\infty} e^{\sqrt{2}x}\frac{t^{3/2}}{\frac{3}{2\sqrt{2}}\log t}U_\phi(t,\sqrt{2}t+x)\le \gamma(r)C(\phi,r).
\end{align}
Letting $r\to\infty$, by the fact that $\lim_{r\to\infty}\gamma(r)=1$, we get that
\begin{align*}
 & \limsup_{r\to\infty}C(\phi,r)\le \liminf_{t\to\infty} e^{\sqrt{2}x}\frac{t^{3/2}}{\frac{3}{2\sqrt{2}}\log t}U_\phi(t,\sqrt{2}t+x)\\
  &\le \limsup_{t\to\infty} e^{\sqrt{2}x}\frac{t^{3/2}}{\frac{3}{2\sqrt{2}}\log t}U_\phi(t,\sqrt{2}t+x)\le \liminf_{r\to\infty}C(\phi,r).
\end{align*}
It follows that $C(\phi):=\lim_{r\to\infty}C(\phi,r)$ exists, and then \eqref{u(t,sqrt(2)t)} follows immediately.
Now we show that $C(\phi)\in(0,\infty)$ if $\phi$ is non-trivial.
In fact, by \eqref{3.6}, we have
$$
0<\gamma(r)^{-1}C(\phi,r)\le C(\phi)\le \gamma(r)C(\phi,r)<\infty.
$$

For any $z$,  it is clear that $U_{\theta_{-z}\phi}(t,x)=U_{\phi}(t,x+z)$, which implies that
$$
C(\phi)e^{-\sqrt{2}z}e^{-\sqrt{2}x}=\lim_{t\to\infty}\frac{t^{3/2}}{\frac{3}{2\sqrt{2}}\log t}U_\phi(t,\sqrt{2}t+x+z)=
C(\theta_{-z}\phi)e^{-\sqrt{2}x},
$$
which further implies that $C(\theta_{-z}\phi)=C(\phi)e^{-\sqrt{2}z}$, that is
\begin{equation}\label{3.3}
  C(\theta_{-z}\phi)=\lim_{r\to\infty}\sqrt{\frac{2}{\pi}}\int_0^\infty U_{\phi}(r,\sqrt{2}r+y+z)ye^{\sqrt{2}y}\,dy=C(\phi)e^{-\sqrt{2}z}.
\end{equation}

(2) Recall that in this part we assume that {\bf (H1)} and {\bf (H3)} hold,
and $\phi\in\mathcal{B}^+(\R)$ satisfies the integrability condition \eqref{initial-cond1'}
at $-\infty$ and that there exists $x_0<0$ such that $\phi$ is bounded on $(-\infty, x_0]$.
Note that , for $t_0>0$,
$U_{\phi}(t+t_0,x+\sqrt{2}t_0)$ and $V_{\phi}(t+t_0,x+\sqrt{2}t_0)$ are the solution to \eqref{KPPpsi} with initial data $U_{\phi}(t_0,x+\sqrt{2}t_0)$ and $V_{\phi}(t_0,x+\sqrt{2}t_0)$ respectively.
By Corollary \ref{initial-cond-V2}, we have that
$U_{\phi}(t_0,x+\sqrt{2}t_0)$ and  $V_{\phi}(t_0,x+\sqrt{2}t_0)$ are bounded
and satisfy
the integrability condition \eqref{initial-cond1} at $\infty$.
It follows from  \eqref{u(t,sqrt(2)t)}  that
\begin{align*}
  \lim_{t\to\infty}\frac{t^{3/2}}{\frac{3}{2\sqrt{2}}\log t}V_{\phi}(t,x+\sqrt{2}t)=\lim_{t\to\infty}\frac{t^{3/2}}{\frac{3}{2\sqrt{2}}\log t}V_{\phi}(t+t_0,x+\sqrt{2}t_0+\sqrt{2}t)=  Ce^{-\sqrt{2}x},
\end{align*}
where
\begin{align*}
  C&=\lim_{r\to\infty}\sqrt{\frac{2}{\pi}}\int_0^\infty V_{\phi}(r+t_0,\sqrt{2}r+\sqrt{2}t_0+y)ye^{\sqrt{2}y}\,dy\\
  &=\lim_{r\to\infty}\sqrt{\frac{2}{\pi}}\int_0^\infty V_{\phi}(r,\sqrt{2}r+y)ye^{\sqrt{2}y}\,dy:=\tl{C}({\phi}).
\end{align*}
Similarly, $$ \lim_{t\to\infty}\frac{t^{3/2}}{\frac{3}{2\sqrt{2}}\log t}U_{\phi}(t,x+\sqrt{2}t)=C(\phi)e^{-\sqrt{2}x}.$$
The proof is now complete.
\hfill
$\Box$

\noindent\textbf{Proof of Theorem \ref{Them:large-der}:}
It is clear that \eqref{lar-der1} follows from \eqref{v(t,sqrt(2)t)} with $\phi=0$.
Now we prove \eqref{lar-der2}. For $t_0>0$, using Proposition \ref{prop:u-psi} with $u(0,x)=V(t_0,\sqrt{2}t_0+x)$,   we get  that
$$
\gamma(r)^{-1}\Psi(r,t,x)\le V(t_0+t,\sqrt{2}t_0+x)\le \gamma(r)\Psi(r,t,x).
$$
By Lemma \ref{initial-cond-V} and the dominated convergence theorem, we have that
\begin{align*}
 &\lim_{t\to\infty}\sqrt{t}e^{(\delta^2/2+\sqrt{2}\delta)t}\Psi(r,t,(\sqrt{2}+\delta)t+x)\\
 &=\sqrt{\frac{2}{\pi}}\delta e^{-(\delta+\sqrt{2})x}e^{-\frac{1}{2}\delta^2 r}\int_0^\infty V(t_0+r,\sqrt{2}r+\sqrt{2}t_0+y)ye^{(\sqrt{2}+\delta)y}\,dy\in(0,\infty).
\end{align*}
Now,
using arguments similar to that used in the proof of Proposition \ref{prop:lim-U} (1), we get that
\begin{align*}
 &\lim_{t\to\infty}\sqrt{t}e^{(\delta^2/2+\sqrt{2}\delta)t}V(t_0+t,(\sqrt{2}+\delta)t+\sqrt{2}t_0+x)\\
&=\sqrt{\frac{2}{\pi}}\delta e^{-(\delta+\sqrt{2})x}\lim_{r\to\infty}e^{-\frac{1}{2}\delta^2 r}\int_0^\infty V(t_0+r,\sqrt{2}r+\sqrt{2}t_0+y)ye^{(\sqrt{2}+\delta)y}\,dy\\
&=\sqrt{\frac{2}{\pi}}\delta e^{\frac{1}{2}\delta^2t_0}e^{-(\delta+\sqrt{2})x}\lim_{r\to\infty}e^{-\frac{1}{2}\delta^2 r}\int_0^\infty V(r,\sqrt{2}r+y)ye^{(\sqrt{2}+\delta)y}\,dy\in(0,\infty),
\end{align*}
where the limit above exists.
Letting $x=\delta t_0$, we get that
\begin{align*}
  &\lim_{t\to\infty}\sqrt{t}e^{(\delta^2/2+\sqrt{2}\delta)t}V(t,(\sqrt{2}+\delta)t)\\
  &=\sqrt{\frac{2}{\pi}}\delta\lim_{r\to\infty}e^{-\frac{1}{2}\delta^2 r}\int_0^\infty V(r,\sqrt{2}r+y)ye^{(\sqrt{2}+\delta)y}\,dy
  :=\hat{C}(\delta).
\end{align*}
It follows that
$$\lim_{t\to\infty}\sqrt{t}e^{(\delta^2/2+\sqrt{2}\delta)t}\P(M_t>(\sqrt{2}+\delta)t)
=\hat{C}(\delta).$$
\hfill
$\Box$

\subsection{The extremal process}

In this subsection
we give the proofs of our main results--Theorems \ref{thrm:tra-wave}, \ref{thrm:tra-wave2} and \ref{them:main3}. Recall that $m(t)$ is defined in \eqref{def-mt}.

\subsubsection{Proofs of Theorems \ref{thrm:tra-wave} and \ref{thrm:tra-wave2}}

\noindent\textbf{Proof of Theorem \ref{thrm:tra-wave}:}
(1)
In this part, we assume that $\phi$ is bounded and satisfies \eqref{initial-cond1'}.
Define
\begin{align}\label{e:Wphi}
w_{\phi}(x):=-\log \E(e^{-C(\phi)Z_\infty e^{-\sqrt{2}x}}).
\end{align}
Recall that (cf. \eqref{Trav-Solution}--\eqref{limit-w})  $w$ is a traveling wave to  \eqref{KPPpsi} and satisfies
\begin{equation}\label{3.8}
  \lim_{x\to\infty}\frac{w_{\phi}(x)}{xe^{-\sqrt{2}x}}=C(\phi).
\end{equation}
Let $\Psi$ be defined by \eqref{e:Psi} with $u$ replaced by $U_\phi$.
We claim that,
for any positive function $z(t)$ with $\lim_{t\to\infty}z(t)=z>0$,
\begin{equation}\label{3.7}
  \lim_{t\to\infty}\Psi(r,t,z(t)+m(t))=C(\phi,r)ze^{-\sqrt{2}z}.
\end{equation}
In fact, for any $t\ge8r$ and $y\ge 0$,
\begin{align}\label{3.10}
  &z(t)^{-1}\frac{t^{3/2}}{\sqrt{2\pi(t-r)}} U_\phi(r,\sqrt{2}r+y)e^{\sqrt{2}y}e^{-\frac{(z(t)+m(t)-\sqrt{2}t-y)^2}{2(t-r)}}\left(1-e^{-2z(t)y/(t-r)}\right)\nonumber\\
  &\le c  U_\phi(r,\sqrt{2}r+y)ye^{\sqrt{2}y}.
\end{align}
Thus, we can apply  the dominated convergence theorem to get that
$$
\lim_{t\to\infty}
z(t)^{-1} e^{\sqrt2 z(t)}\Psi(r,t,z(t)+m(t))=C(\phi,r),
$$
which is the same as \eqref{3.7}.
Put
$$f^{(r,t)}(x)=|\Psi(r,t,x+m(t))-C(\phi,r)xe^{-\sqrt{2}x}|,\quad x>0.$$
 It follows from the key inequality \eqref{bound-u3}
that $U_\phi(t,x)\le (\|\phi\|_\infty\vee 1).$
Applying Proposition \ref{prop:u-psi},
we get that, for $r$ large enough, $x>9r$ and $t\ge 8r$,
$$
U_\phi(t,m(t)+x)\le \gamma(r)\Psi(r,t,x+m(t)).
$$
Thus, for any $t\ge 8r$,
\begin{align*}
 U_\phi(t,m(t)+x)\le
  \gamma(r)C(\phi,r)xe^{-\sqrt{2}x}1_{x>9r}
 +\gamma(r)f^{(r,t)}(x)1_{x>9r}+(\|\phi\|_\infty\vee 1)1_{x\le 9r}.
\end{align*}
Let
$u_{1,r}(s,x)$, $u_{2,r}(s,x)$ and
$v_{r,t}(s,x) $ be the solutions to the KPP equation \eqref{KPPpsi} with initial conditions
$C(\phi)xe^{-\sqrt{2}x}1_{x>9r}$,
$(\|\phi\|_\infty\vee 1)1_{x\le 9r}$ and $\gamma(r)f^{(r,t)}(x)1_{x>9r}$ respectively.
Then, by Lemma \ref{lem:compare}, we have
$$U_\phi(t+s,m(t)+x)\le \left(\frac{\gamma(r)C(\phi,r)}{C(\phi)}\vee1\right)
u_{1,r}(s,x)+u_{2,r}(s,x)+v_{r,t}(s,x).$$
Let $a(r):=\frac{\gamma(r)C(\phi,r)}{C(\phi)}\vee1$. Applying \eqref{F-K-u} with $r=0$ and using the fact that $k(\lambda)\le 1$,  we get that
$$v_{r,t}(s,x)\le e^s\gamma(r)\mE_x(f^{(r,t)}(B_s)1_{B_s>9r}).$$
Thus,
\begin{align*}
& U_\phi(t+s,m(t+s)+x(t+s))\\
\le& a(r)u_{1,r}(s,x(t+s)+m(t+s)-m(t))
+u_{2,r}(s,x(t+s)+m(t+s)-m(t))\\
 &+e^s\gamma(r)\mE(f^{(r,t)}(m(t+s)-m(t)+x(t+s)+B_s),m(t+s)-m(t)+x(t+s)+B_s>9r).
\end{align*}
Letting $t\to\infty$ and using \eqref{3.7}, we get
$$\limsup_{t\to\infty}U_\phi(t,m(t)+x(t))\le a(r)u_{1,r}(s,x+\sqrt{2}s)+u_{2,r}(s,x+\sqrt{2}s).$$
Since $(\|\phi\|_\infty\vee 1)1_{x<9r}$ satisfies \eqref{initial-cond1}, we have by Proposition \ref{prop:lim-U}.(1) that $u_{2,r}(s,x+\sqrt{2}s)\to 0$ as $s\to\infty$.  Since $C(\phi)xe^{-\sqrt{2}x}/w_{\phi}(x)\to 1$ as $x\to\infty$, by \cite[Lemma 3.4]{Bramson}
(Note that our \eqref{limit-w}, which is exactly \cite[(1.13)]{Bramson}, holds under {\bf (H1)-(H2)}, thus
\cite[Lemma 3.4]{Bramson} holds under {\bf (H1)-(H2)}),
we get that
$$u_{1,r}(s,x+\sqrt{2}s)\to w_\phi(x),\quad s\to\infty.$$
Now letting $s\to\infty$ and then $r\to\infty$, we get that
 $$\limsup_{t\to\infty}U_\phi(t,m(t)+x(t))\le  \lim_{r\to\infty}a(r)w_\phi(x)=w_\phi(x).$$
On the other hand,
 $$\gamma(r)^{-1}C(\phi,r)xe^{-\sqrt{2}x}1_{x>9r}\le U_\phi(t,m(t)+x)+\gamma(r)^{-1}f^{(r,t)}(x)1_{x>9r}.$$
Using arguments similar as above, we can get that
$$\liminf_{t\to\infty}U_\phi(t,m(t)+x(t))\ge  w_\phi(x).$$
Therefore, we have
\begin{equation}\label{3.9}
  \lim_{t\to\infty}U_\phi(t,m(t)+x(t))=  w_\phi(x).
\end{equation}

(2)
Recall that, in this part, $\phi$ is not necessarily  bounded.
Applying \eqref{3.9} to $(t,x)\to
V_\phi(t+t_0,\sqrt{2}t_0+x)$
 and Proposition \ref{prop:lim-U} (2), we get that
$$\lim_{t\to\infty}V_\phi(t+t_0,m(t)+\sqrt{2}t_0+x(t))=-\log \E(e^{-\tl{C}(\phi)Z_\infty e^{-\sqrt{2}x}}).$$
Since $m(t+t_0)-m(t)-\sqrt{2}t_0+x(t)\to x$,  we get that
$$\lim_{t\to\infty}V_\phi(t+t_0,m(t+t_0)+x(t))=-\log \E(e^{-\tl{C}(\phi)Z_\infty e^{-\sqrt{2}x}}),$$
which implies the desired result.
Similarly, applying Corollary \ref{initial-cond-V2} and \eqref{3.9} to $(t,x)\to
U_\phi(t+t_0,\sqrt{2}t_0+x)$, it is clear that \eqref{3.9} also holds for $\phi\in\mathcal{B}^+(\R)$.
The proof is now complete.
\hfill$\Box$

Using Theorem \ref{thrm:tra-wave}, we get the convergence of the Laplace transforms. To obtain weak convergence, we need to show the continuity of $C(\phi)$ and $\tl C(\phi)$.

\begin{lemma}\label{cont-C}
Assume that $\bf{(H1)}$ and $\bf{(H2)}$ hold. Then for any $\phi\in\mathcal{H}$,
\begin{equation}\label{cont-C1}
\lim_{\lambda \downarrow0}C(\lambda\phi)=C(0)=0.
\end{equation}
If, in addition, $\bf{(H3)}$ holds, then for any $\phi\in\mathcal{H}$,
\begin{equation}\label{cont-C2}
  \lim_{\lambda\downarrow0}\tl C(\lambda\phi)=\tl C_0.
\end{equation}
\end{lemma}
\textbf{Proof:}
For any $\phi\in\mathcal{H}$, choose $m_\phi$ such that $\phi(x)=0$ for all $x<m_\phi$. Then we have for all $N>0$,
$$\E(\exp\{-\lambda\langle\phi,\cE_t\rangle\})\ge \E(\exp\{-\lambda\|\phi\|_\infty\cE_t(m_\phi,\infty)\})\ge e^{-\lambda\|\phi\|_\infty N}\P(\cE_t(m_\phi,\infty)\le N).$$
Letting $t\to\infty$, $\lambda\to0$ and then $N\to\infty$,
by Theorem \ref{thrm:tra-wave} we see that,
to prove \eqref{cont-C1}, it suffices to show that
\begin{equation}\label{tight}
  \lim_{N\to\infty}\limsup_{t\to\infty}\P(\cE_t(m_\phi,\infty)> N)=0.
\end{equation}
Let $g(x)={\bf 1}_{(0,\infty)}(x)$, then $u_g(t,x)=-\log\E(\exp\{-X_t(-x,\infty)\})$ is increasing on $x$. For any $n\ge1$,
\begin{align*}
  &\E(\cE_t(m_\phi,\infty)> N,\exp\{-\langle \theta_{-n}g,\cE_{t+1}\})\\
 & =\E\left(\cE_t(m_\phi,\infty)> N,\E_{X_t}\Big(\exp\Big\{-\int g(y-n-m(t+1))X_{1}(dy)\Big\}\Big)\right) \\
   &= \E\left(\cE_t(m_\phi,\infty)> N,\exp\Big\{-\int u_g(1,x-n-m(t+1)+m(t))\cE_{t}(dx)\Big\}\right)\\
   &\le \E\left(\cE_t(m_\phi,\infty)> N,\exp\Big\{-u_g(1,m_{\phi}-n-m(t+1)+m(t))\cE_t(m_\phi,\infty)\Big\}\right)\\
   &\le \exp\{-u_g(1,m_{\phi}-n-m(t+1)+m(t))N\}.
\end{align*}
Thus, we get that
\begin{align*}
  &\limsup_{t\to\infty}\P(\cE_t(m_\phi,\infty)> N)\\
    &\le  \limsup_{t\to\infty}\E(\cE_t(m_\phi,\infty)> N,\exp\{-\langle \theta_{-n}g,\cE_{t+1}\})\\
&\quad +1-\lim_{t\to\infty}\E(\exp\{-\langle \theta_{-n}g,\cE_{t+1}\})\\
   &\le \exp\{-u_g(1,m_{\phi}-n-\sqrt{2})N\}+1-\E(\exp\{-C(g)e^{-\sqrt{2}n}Z_\infty\}).
\end{align*}
Letting $N\to\infty$ and then $n\to\infty$, \eqref{tight} follows immediately.
Thus \eqref{cont-C1} is valid.

Now we prove \eqref{cont-C2} under the additional assumption $\bf{(H3)}$.
It is clear that
$$0\le \P(M_t-m(t)\le 0)-\E(\exp\{-\lambda\langle\phi,\cE_t\rangle\},M_t-m(t)\le 0)\le 1-\E(\exp\{-\lambda\langle\phi,\cE_t\rangle\}).$$
Thus, by \eqref{cont-C1} and
Theorem \ref{thrm:tra-wave}.(2) with $x(t)=0$, we get that
$$\E(\exp\{-\lim_{\lambda\to0}\tl C(\lambda\phi)Z_\infty\})=\E(\exp\{-\tl C_0Z_\infty\}).$$
Now \eqref{cont-C2} follows immediately.
The proof is now complete.
\hfill$\Box$

For any $t>0$, we define $\bar{\cal E}_t:={\cal T}_{-\sqrt{2}t}X_t$.

\begin{proposition}\label{prop-lim2}
Assume that $\bf{(H1)}$ and $\bf{(H3)}$ hold. For any $z\in \R$,
 conditioned on $\{M_t>\sqrt{2}t+z\}$, $(\bar{\cal E}_t-z, M_t-\sqrt{2}t-z)$ converges in distribution to a limit (independent of $z$) $(\bar{\cal E}_\infty, Y)$, where $Y$ is an exponential random variable with parameter $\sqrt{2}$ and for any $\phi\in{\cal C}_c^+(\R)$,
$$
\E\Big[\exp\Big\{-\int_{\R}\phi(y)\bar{\cal E}_\infty(dy)\Big\},Y>x\Big]=\frac{\tl{C}(\theta_x\phi)e^{-\sqrt{2}x}-C(\phi)}{\tl{C}_0}, \quad x>0.
$$
Moreover, $$(X_t-M_t,M_t-\sqrt{2}t-z)|_{M_t>\sqrt{2}t+z}$$
converges in law to $(\Delta,Y)$, where the random measure $\Delta=\bar{\cal E}_\infty-Y$ is independent of $Y$.
\end{proposition}

\begin{remark}
Define $Y_t:=X_t(\sqrt{2 }t,\infty)=\langle h,\bar{\cE}_t\rangle$, where $h(x)=1_{(0,\infty)}(x)$.
It follows from Proposition \ref{prop-lim2} that $Y_t|_{Y_t>0}$ converges weakly.
\end{remark}

\noindent\textbf{Proof of Proposition \ref{prop-lim2}:}
First, we show that $M_t-\sqrt{2}t-z|_{M_t>\sqrt{2}t+z}$ converges in distribution to an exponential random variable with parameter $\sqrt{2}$.
In fact,
by \eqref{v(t,sqrt(2)t)} with $\phi=0$, we get that for any $x>0$,
\begin{align*}
  \lim_{t\to\infty}\P(M_t-\sqrt{2}t-z>x|M_t>\sqrt{2}t+z) & =\lim_{t\to\infty}\frac{\frac{t^{3/2}}{\frac{3}{2\sqrt{2}}\log t}\P(M_t>\sqrt{2}t+z+x)}{\frac{t^{3/2}}{\frac{3}{2\sqrt{2}}\log t}\P(M_t>\sqrt{2}t+z)}=e^{-\sqrt{2}x}.
\end{align*}
For any $x>0$ and $\phi\in\mathcal{B}_b^+(\R)$ satisfying
the integrability condition \eqref{initial-cond1'} at $-\infty$,
  applying Proposition \ref{prop:lim-U} several times, we get
\begin{align*}
  &\lim_{t\to\infty}\E\Big(e^{-\int_{\R} \phi(y-\sqrt{2}t-z)\,X_t(dy)},M_t> \sqrt{2}t+z+x |M_t>\sqrt{2}t+z\Big)\\
  &=\lim_{t\to\infty}\frac{\E\Big(e^{-\int_{\R } \phi(y-\sqrt{2}t-z)\,X_t(dy)},M_t>\sqrt{2}t+z+x\Big)}{\P(M_t>\sqrt{2}t+z)}  \\
  &=\lim_{t\to\infty}\frac{1-\E(e^{-\int_{\R} \phi(y-\sqrt{2}t-z)\,X_t(dy)},M_t
  \le \sqrt{2}t+z+x)}{\P(M_t>\sqrt{2}t+z)}-\lim_{t\to\infty}\frac{1-\E(e^{-\int_{\R} \phi(y-\sqrt{2}t-z)\,X_t(dy)})}{\P(M_t>\sqrt{2}t+z)}\\
  &=\lim_{t\to\infty}\frac{V_{\theta_x\phi}(t,\sqrt{2}t+z+x)}{V(t,\sqrt{2}t+z)}-\lim_{t\to\infty}\frac{U_{\phi}(t,\sqrt{2}t+z)}{V(t,\sqrt{2}t+z)}\\
  &=\frac{\tl{C}(\theta_x\phi)e^{-\sqrt{2}x}-C(\phi)}{\tl{C}_0},
\end{align*}
where in the second to last equality above, we used L'Hospital's rule and the facts
that $V_{\theta_x\phi}(t,\sqrt{2}t+z+x)\to 0$ and $U_{\phi}(t,\sqrt{2}t+z)\to 0$
(which are consequences of \eqref{v(t,sqrt(2)t)} and \eqref{u(t,sqrt(2)t)} respectively).
Now applying Lemma \ref{cont-C},  we get that for any
$\phi\in\mathcal{C}_c^+(\R)$,
$$(\langle\phi,\bar{\cE}_t-z\rangle,M_t-\sqrt{2}t-z)|_{M_t>\sqrt{2}t+z}$$
converges jointly in distribution. Thus the limit has the form $(\langle\phi,\bar{\cE}_\infty\rangle,Y)$ (independent of $z$), where the random measure $\bar{\cE}_\infty\in\cM_R(\R)$.

It follows by \cite[Lemma 4.13]{ABK} that, conditioned on $M_t>\sqrt{2}t+z$, $X_t-M_t$ converges in law to $\bar{\cal E}_\infty-Y$. Thus,
\begin{align*}
  &\E\Big(e^{-\int_{\R} \phi(y-M_t)\,X_t(dy)},M_t-\sqrt{2}t-z> x |M_t>\sqrt{2}t+z\Big)\\
   &= \E\Big(e^{-\int_{\R} \phi(y-M_t)\,X_t(dy)}|M_t>\sqrt{2}t+z+x\Big)\P(M_t>\sqrt{2}t+z+x|M_t>\sqrt{2}t+z) \\
   & \to \E(e^{-\int_{\R} \phi(y-Y)\,\bar{\cal E}_\infty(dy)})\P(Y>x).
\end{align*}
The desired independence result follows immediately.
\hfill$\Box$

\noindent\textbf{Proof of Theorem \ref{thrm:tra-wave2}:}
The weak convergence of $\cE_t$ and \eqref{laplace-E-M} follow immediately from Theorem \ref{thrm:tra-wave} (1) and Lemma \ref{cont-C}.  Now we assume
that $\bf{(H3)}$ also holds and prove the second assertion of Theorem \ref{thrm:tra-wave2}.
For any $\phi\in \mathcal{C}_c^+(\R)$,
choose $m_\phi$ such that $\phi(y)=0$ for all $y<m_\phi$.
Then we have
\begin{align}
  \frac{C(\phi)}{\tl{C}_0}&=\lim_{t\to\infty}\frac{1-\E \Big[e^{-\int \phi(y-\sqrt{2}t)\,X_t(dy)}\Big]}{\P(M_t>\sqrt{2}t)}\nonumber\\
  &=\lim_{t\to\infty}\frac{\E\Big[1- e^{-\int \phi(y-\sqrt{2}t)\,X_t(dy)},M_t>\sqrt{2}t+m_{\phi}\Big]}{\P(M_t>\sqrt{2}t)}\nonumber\\
  &=\lim_{t\to\infty}\E\Big[1- e^{-\int \phi(y-\sqrt{2}t)\,X_t(dy)}|M_t>\sqrt{2}t+m_{\phi}\Big]\frac{\P(M_t>\sqrt{2}t+m_{\phi})}{\P(M_t>\sqrt{2}t)}\nonumber\\
  &=e^{-\sqrt{2}m_{\phi}}\lim_{t\to\infty}\E\Big[1- e^{-\int \phi(y+m_{\phi}-\sqrt{2}t-m_{\phi})\,X_t(dy)}|M_t>\sqrt{2}t+m_{\phi}\Big]\nonumber\\
  &=e^{-\sqrt{2}m_{\phi}}\E\Big[1- e^{-\int \phi(y+m_{\phi})\,\bar{\mathcal{E}}_\infty(dy)}\Big]
  =e^{-\sqrt{2}m_{\phi}}\int_{0}^\infty \sqrt{2}e^{-\sqrt{2}x}\E\Big[1- e^{-\int \phi(y+m_\phi+x)\,\Delta(dy)}\Big]\,dx\nonumber\\
  &=\int_{-\infty}^\infty \sqrt{2}e^{-\sqrt{2}x}\E\Big[1- e^{-\int \phi(y+x)\,\Delta(dy)}\Big]\,dx,\label{c-phi-c-0}
\end{align}
where in the first, fifth and sixth equality we used Proposition \ref{prop-lim2}, and in the fourth
equality we used Proposition  \ref{prop:lim-U}.
By the definition of $\sum_j\Delta_j(dx+e_j)$, we deduce that
\begin{align}
  &\E\Big(e^{-\sum_{j}\phi(y+e_j)\Delta_j(dy)}\Big) =\E\prod_j\Big[\E\Big(e^{-\int \phi(y+x)\Delta(dy)}\Big)\Big]_{x=e_j} \nonumber\\
   & =\exp\left\{-\int \Big(1-\E\Big(e^{-\int\phi(y+x)\Delta(dy)}\Big)\Big)\tl{C}_0Z_\infty \sqrt{2}e^{-\sqrt{2}x}\,dx\right\}=\exp\{-C(\phi)Z_\infty\}.\label{laplace-delta-e}
\end{align}
The proof is now complete.
\hfill$\Box$

\subsubsection{Proof of Theorem \ref{them:main3}}

\begin{lemma}\label{lem:main3}
Assume that {\bf (H1)} and  {\bf (H3)} hold, and that
$\phi\in{\cal B}^+(\R)$ satisfies the integrability condition \eqref{initial-cond1'}
at $-\infty$ and that there exists $x_0<0$ such that $\phi$ is bounded on $(-\infty, x_0]$.
Then
\begin{align*}
  \lim_{s\to\infty}\lim_{t\to\infty}\E\Big[e^{-\int \phi(y-m(t)-\frac{1}{\sqrt{2}}\log Z_s)X_t(dy)}e^{-\theta Z_s},Z_s>0\Big]
    =e^{-C(\phi)}\E\Big[e^{-\theta Z_\infty},Z_\infty>0\Big].
\end{align*}
\end{lemma}
\textbf{Proof:}
By the Markov property, we have for $s<t$,
\begin{align*}
  &\E\Big[\exp\Big\{-\int \phi(y-m(t)-\frac{1}{\sqrt{2}}\log Z_s)X_t(dy)\Big\}\exp\{-\theta Z_s\},Z_s>0\Big]\\
  &=\E\Big[\exp\Big\{-\int U_\phi(t-s,m(t)+\frac{1}{\sqrt{2}}\log Z_s-y)X_s(dy)\Big\}\exp\Big\{-\theta Z_s\Big\},Z_s>0\Big].
\end{align*}
Now applying Theorem \ref{thrm:tra-wave} (2)
and \eqref{e:Wphi}, we get that as $t\to\infty$,
\begin{align*}
&\E\Big[\exp\Big\{-\int \phi(y-m(t)-\frac{1}{\sqrt{2}}\log Z_s)X_t(dy)\Big\}\exp\{-\theta Z_s\},Z_s>0\Big]\\
&\to \E\Big[\exp\Big\{-\int w_\phi(\sqrt{2}s+\frac{1}{\sqrt{2}}\log Z_s-y)X_s(dy)\Big\}\exp\Big\{-\theta Z_s\Big\},Z_s>0\Big].
\end{align*}
For any $L>0$, define $A(s,L):=\{Z_s>0,\log Z_s\in[-L,L], M_s\le \sqrt{2}s-\log s\}$. Then
\begin{align*}
  &\E\Big[\exp\Big\{-\int w_\phi(\sqrt{2}s+\frac{1}{\sqrt{2}}\log Z_s-y)X_s(dy)\Big\}\exp\Big\{-\theta Z_s\Big\},Z_s>0\Big]  \\
  & \le \E\Big[\exp\Big\{-\int w_\phi(\sqrt{2}s+\frac{1}{\sqrt{2}}\log Z_s-y)X_s(dy)\Big\}\exp\Big\{-\theta Z_s\Big\},A(s,L)\Big]\\
  &\quad +\P(Z_s>0,|\log Z_s|>L)+\P(M_s>\sqrt{2}s-\log s):=(I)+(II)+(III).
\end{align*}
Since $\{Z_\infty=0\}=\cS^c$ and $\{Z_\infty>0\}=\cS$ a.s.,  we have
\begin{align}\label{e:II}
&\lim_{L\to\infty}\limsup_{s\to\infty}\P(Z_s>0,|\log Z_s|>L)\nonumber\\
&\le \lim_{L\to\infty}\limsup_{s\to\infty}\Big(\P(Z_s>0,|\log Z_s|>L,\cS)+\P(X_s\neq 0,\cS^c)\Big)\nonumber\\
&\le \lim_{L\to\infty}\P(|\log Z_\infty|\ge L,\cS)=0.
\end{align}
By \eqref{1.1}, we have
\begin{align}\label{e:III}
\lim_{s\to\infty}\P(M_s>\sqrt{2}s-\log s)=0.
\end{align}
Now we consider $(I)$.
Since $\frac{w_\phi(x)}{xe^{-\sqrt{2}x}}\to C(\phi)$, as $x\to\infty$, and on
$A(s,L)$, for $y\in \mbox{supp}\ X_s$,
$\sqrt{2}s+\frac{1}{\sqrt{2}}\log Z_s-y\ge \log s-L/\sqrt{2}$, thus for any $\epsilon>0,$ there exists $N$
such that for $s>N$,
\begin{align*}
  &(1-\epsilon)C(\phi)\int(\sqrt{2}s+\frac{1}{\sqrt{2}}\log Z_s-y)e^{-\sqrt{2}(\sqrt{2}s+\frac{1}{\sqrt{2}}\log Z_s-y)}X_s(dy)\\
  &\le \int w_\phi(\sqrt{2}s+\frac{1}{\sqrt{2}}\log Z_s-y)X_s(dy)\\
  &\le (1+\epsilon)C(\phi)\int (\sqrt{2}s+\frac{1}{\sqrt{2}}\log Z_s-y)e^{-\sqrt{2}(\sqrt{2}s+\frac{1}{\sqrt{2}}\log Z_s-y)}X_s(dy).
 \end{align*}
 Note that on  $A(s,L)$, for $s$ large enough, $\frac{|\log Z_s|}{\sqrt{2}(\sqrt{2}s-y)}\le \frac{L}{\sqrt{2}\log s}\le \epsilon.$
 Thus $(I)$ is less than or equal to
 \begin{align}\label{e:Iu}
  &\E\Big[\exp\Big\{-(1-\epsilon)^2C(\phi)(Z_s)^{-1}\int (\sqrt{2}s-y)e^{-\sqrt{2}(\sqrt{2}s-y)}X_s(dy)\Big\}\exp\Big\{-\theta Z_s\Big\},A(s,L)\Big]\nonumber\\
   &\le \exp\Big\{-(1-\epsilon)^2C(\phi)\Big\}\E\Big[\exp\Big\{-\theta Z_s\Big\},Z_s>0\Big].
 \end{align}
 Similarly,
 \begin{align}\label{e:Il}
   &(I)\ge \exp\Big\{-(1+\epsilon)^2C(\phi)\Big\}\E\Big[\exp\Big\{-\theta Z_s\Big\},A(s,L)\Big].
 \end{align}
Combining \eqref{e:II}--\eqref{e:Il},
letting $s\to\infty$, then $L\to\infty,$ and then
$\epsilon\to0$, we get that
 \begin{align*}
   &\lim_{s\to\infty}\E\Big[\exp\Big\{-\int w_\phi(\sqrt{2}s+\frac{1}{\sqrt{2}}\log Z_s-y)X_s(dy)\Big\}\exp\Big\{-\theta Z_s\Big\},Z_s>0\Big] \\
   & =\exp\Big\{-C(\phi)\Big\}\E\Big[\exp\Big\{-\theta Z_\infty\Big\},Z_\infty>0\Big].
 \end{align*}
The proof is now complete.
\hfill
$\Box$

\noindent\textbf{Proof of Theorem \ref{them:main3}:}
Using  arguments similar to that leading to \eqref{laplace-delta-e},
we get that for any $\phi\in{\cal C}_c^+(\R)$,
\begin{align*}
  &\E\Big(e^{-\sum_{j}\phi(y+e_j')\Delta_j(dy)}\Big)
   =\exp\left\{-\int \Big(1-\E\Big(e^{-\int\phi(y+x)\Delta(dy)}\Big)\Big)\tl{C}_0 \sqrt{2}e^{-\sqrt{2}x}\,dx\right\}=\exp\{-C(\phi)\}.
\end{align*}
Since $Z_t\to Z_\infty$, we only need to prove that, for any  $\phi\in{\cal C}_c^+(\R)$ and $\theta\ge0$,
\begin{align}
  &\lim_{t\to\infty}\E\Big[e^{-\int \phi(y)\cE^*_t(dy)}e^{-\theta Z_\infty},Z_\infty>0\Big]
=\exp\{-C(\phi)\}\E[\exp\{-\theta Z_\infty\},Z_\infty>0].\label{4.2}
\end{align}

{\bf Step 1}
Define, for any $b>1$,
$$g_b(x):=\left\{
            \begin{array}{ll}
              0, & \hbox{$|x|>b$;} \\
              1, & \hbox{$|x|<b-1$;} \\
              \hbox{linear}, & \hbox{otherwise.}
            \end{array}
          \right.
$$
It is clear that $|g_b(x)-g_b(y)|\le |x-y|$.

First,  we consider
$\phi(x)=|f(x)|g_b(x)$ where
$f(x)=\sum_{i=1}^n \theta_i e^{\beta_ix}$ and $\theta_i,\beta_i\in\R.$
Let $\overline{f}(x):=\sum_{i=1}^n |\theta_i| e^{\beta_ix}$.
It is clear that $\phi\in{\cal C}_c^+(\R)$.
By Lemma \ref{lem:main3}, to prove \eqref{4.2}, it suffices to show that
\begin{align}\label{4.1}
  \lim_{s\to\infty}\limsup_{t\to\infty}\Big|&\E\Big[e^{-\int \phi(y-m(t)-\frac{1}{\sqrt{2}}\log Z_\infty)X_t(dy)}e^{-\theta Z_\infty},Z_\infty>0\Big]\nonumber\\
&-\E\Big[e^{-\int \phi(y-m(t)-\frac{1}{\sqrt{2}}\log Z_s)X_t(dy)}e^{-\theta Z_s},Z_s>0\Big]\Big|=0.
\end{align}
For any  $K>0$ and $M>0$, let $\overline{f}^M(y)=\overline{f}(y){\bf 1}_{|y|\le M+b}$ and
$$
A(s,t,K,M)=\{\langle \overline{f}^M,\cE_t\rangle> K\}\cup\left\{\frac{1}{\sqrt{2}}|\log Z_\infty|>M\right\}\cup\left\{\frac{1}{\sqrt{2}} |\log Z_s|>M\right\}.
$$
Since $\langle \overline{f}^M,\cE_t\rangle$ converges weakly and $Z_s\to Z_\infty$, a.s.,
for any $\epsilon>0$, there exist $K,M$ such that
\begin{equation}\label{3.5}
    \lim_{s\to\infty}\limsup_{t\to\infty}\P(A(s,t,K,M),Z_\infty>0)<\epsilon.
\end{equation}
Note that, if  $|y|>|x_1|\vee |x_2|+b$, $\phi(y-x_1)-\phi(y-x_2)=0$; otherwise,
\begin{align*}
& |\phi(y-x_1)-\phi(y-x_2)|\\
&\le |f(y-x_1)-f(y-x_2)|+|f(y-x_2)||g_b(y-x_1)-g_b(y-x_2)|\\
&\le \overline{f}(y)\left[\sum_{j}|e^{-\beta_jx_1}-e^{-\beta_jx_2}|+\sum_{j}e^{-\beta_jx_2}|x_1-x_2|\right]=:\overline{f}(y)H(x_1,x_2).
\end{align*}
By the inequality $|e^{-x}-e^{-y}|\le 1-e^{-|x-y|}$ for any $x,y>0$, we get that on $A(s,t,K,M)^c\cap\{Z_s>0,Z_\infty>0\}$,
\begin{align*}
  &\Big|e^{-\int \phi(y-m(t)-\frac{1}{\sqrt{2}}\log Z_\infty)X_t(dy)}e^{-\theta Z_\infty}-e^{-\int \phi(y-m(t)-\frac{1}{\sqrt{2}}\log Z_s)X_t(dy)}e^{-\theta Z_s}\Big|\\
  & \le 1-\exp\Big\{-\theta|Z_s-Z_\infty|-\langle
  \overline{f}^M, \cE_t\rangle H\left(\frac{1}{\sqrt{2}}\log Z_\infty,\frac{1}{\sqrt{2}}\log Z_s\right)\Big\}\\
&\le 1-\exp\Big\{-\theta|Z_s-Z_\infty|-KH\left(\frac{1}{\sqrt{2}}\log Z_\infty,\frac{1}{\sqrt{2}}\log Z_s\right)\Big\}.
\end{align*}
Since $Z_s\to Z_\infty$, the left hand side of \eqref{4.1} is no more than
\begin{align*}
   &\lim_{s\to\infty}\limsup_{t\to\infty}\Big[\P(Z_s\le 0,Z_\infty>0)+\P(A(s,t,K,M),Z_\infty>0)\\
 &+\E
\left(1-\exp\Big\{-\theta|Z_s-Z_\infty|-KH\left(\frac{1}{\sqrt{2}}\log Z_\infty,\frac{1}{\sqrt{2}}\log Z_s\right)\Big\},Z_s>0,Z_\infty>0\right)\Big]\\
&\le \epsilon.
\end{align*}
Now \eqref{4.1} follows immediately. Thus, the result holds for
$\phi(x)$ of the form specified at the beginning of this step.

{\bf Step 2} We will show that \eqref{4.2} holds for $\phi\in {\cal C}_c^+(\R)$. Choose $b>1$  such that $\phi(x)=0$ for $|x|>b-1$. According to the Stone-Weierstrass theorem, for any $n\ge 1$, there exists a polynomial  $Q_{n,b}$ such that
$$\sup_{y\in[e^{-b},e^{b}]}|Q_{n,b}(y)-\phi(\log y)|\le n^{-1},$$
which is equivalent to that
$$\sup_{y\in[-b,b]}|Q_{n,b}(e^{y})-\phi(y)|\le n^{-1}.$$
Let
$\phi_{n,b}(y):=|Q_{n,b}(e^{y})|g_b(y)$,
it is clear that all the functions $\phi_{n, b}$ satisfy the conditions in Step 1, and $|\phi_{n,b}(y)-\phi(y)|\le n^{-1}g_b(y).$ Thus
\begin{align*}
  &\Big|\E\Big[e^{-\int \phi(y)\cE^*_t(dy)}e^{-\theta Z_\infty},Z_\infty>0\Big]-\E\Big[e^{-\int \phi_{n,b}(y)\cE^*_t(dy)}e^{-\theta Z_\infty},Z_\infty>0\Big]\Big| \\
  & \le \E\Big[1-e^{-\int |\phi(y)-\phi_{n,b}(y)|\cE^*_t(dy)},Z_\infty>0\Big]\\
&\le \E\Big[1-e^{-n^{-1}\int g_b(y)\cE^*_t(dy)}, Z_\infty>0\Big].
\end{align*}
In Step 1, we have shown that,
$$\lim_{n\to\infty}\lim_{t\to\infty}\E\Big[1-e^{-n^{-1}\int g_b(y)\cE^*_t(dy)}, Z_\infty>0\Big]=\lim_{n\to\infty}(1-\exp\{-C(n^{-1}g_b)\})\P(Z_\infty>0)=0.$$
Thus we have
\begin{align*}
  &\lim_{t\to\infty}\E\Big[e^{-\int \phi(y)\cE^*_t(dy)}e^{-\theta Z_\infty},Z_\infty>0\Big]
  =\lim_{n\to\infty}\lim_{t\to\infty}\E\Big[e^{-\int \phi_{n,b}(y)\cE^*_t(dy)}e^{-\theta Z_\infty},Z_\infty>0\Big]\\
&=\lim_{n\to\infty}\exp\{-C(\phi_{n,b})\}\E\Big[e^{-\theta Z_\infty},Z_\infty>0\Big].
\end{align*}
Since $|\phi_{n,b}(y)-\phi(y)|\le n^{-1}g_b(y)$, by Lemmas \ref{lem:compare} and  \ref{cont-C},
we have
 $$|C(\phi_{n,b})-C(\phi)|\le C(n^{-1}g_b)\to 0,\quad n\to\infty.$$
Thus, \eqref{4.2} is valid for all $\phi\in\mathcal{C}_c^+(\R)$.
\hfill
$\Box$

In fact, it is not surprising that $\cE^*_\infty$ has the decomposition \eqref{decom:E*}.
A random measure $M$ is said to be exp-$\sqrt{2}$-stable if any $a,b$ satisfying $e^{\sqrt{2}a}+e^{\sqrt{2}b}=1$, it holds that
$${\cal T}_{a}M+{\cal T}_{b}\hat{M}\overset{d}{=}M,$$
where $\hat{M}$ is an independent copy of $M$.
The following proposition shows that  $\cE^*_\infty$
satisfies the  exp-$\sqrt{2}$-stability:
\begin{proposition}\label{prop-exp-stab}
Under $\bf{(H1)}$ and $\bf{(H3)}$, $\cE^*_\infty$
satisfies the  exp-$\sqrt{2}$-stability.
\end{proposition}
{\bf Proof:} The Laplace transform of ${\cal T}_a\cE^*_\infty$ is given by
$$
\E(\exp\{-\langle \phi, {\cal T}_{a}\cE^*_\infty\})
=\E(\exp\{-\langle \theta_a\phi, \cE^*_\infty\})=\exp\{-C(\theta_a\phi)\}=
\exp\{-C(\phi)e^{\sqrt{2}a}\}, \quad \phi\in \mathcal C_c^+(\R).
$$
Therefore, the desired result follows.
\hfill$\Box$

\begin{remark}
Let $M_1,\dots, M_n$ be a sequence of i.i.d. random measures with the same law as
${\cal T}_{-\log n/\sqrt{2}}\cE^*_\infty.$
Then, by Proposition \ref{prop-exp-stab}, $\cE^*_\infty$ is equal in law to $M_1+\cdots+M_n$.
Thus $\cE^*_\infty$ is infinitely divisible. Applying  \cite[Theorem 3.1]{Maillard}, we get that for any $\phi\in\mathcal{C}_c^+(\R)$,
\begin{align*}
  C(\phi)=-\log \E[\exp\{-\langle \phi,\cE^*_\infty\rangle\}]=&c\int_{\R} \phi(x)e^{-\sqrt{2}x}\,dx\\
 & +\int_\R e^{-\sqrt{2}x}\int_{\cM_R(\R)\setminus \{0\}}[1-\exp\{-\langle \phi,\mu\rangle\}]\,{\cal T}_x\Lambda(d\mu)\,dx,
\end{align*}
for some constant $c>0$ and some measure $\Lambda$ on $\cM_R(\R)\setminus \{0\}$ with the property that for every bounded Borel set $A\subset \R$,
$$
\int_\R e^{-\sqrt{2}x}\int_{\cM_R(\R)\setminus\{0\}}[1\wedge \mu(A-x)]\Lambda(d\mu)\,dx<\infty.
$$
Now we choose a function  $\phi\in\mathcal{C}_c^+(\R)$ such that
$\phi(x)=0$ for any $x<0$.
It is clear that $U_{\lambda\phi}(t,x)\le V(t,x)$. Under {\bf (H1)} and {\bf (H3)}, it holds that $C(\lambda\phi)\le \tl{C}_0\in(0,\infty)$ for any $\lambda>0$.
This implies that $c=0$. Thus
$$\cE^*_\infty\overset{d}{=}\sum_{j}{\cal T}_{\xi_j}D_j,$$
where $\sum_j\delta_{(x_j,D_j)}$ is Poisson point process with intensity measure $e^{-\sqrt{2}x}\,dx\,\Lambda(d\mu).$
Theorem \ref{them:main3} says that
$\Lambda(d\mu)=\sqrt{2}\tl{C}_0\P(\Delta\in d\mu)$
where $\Delta$ is
the limit of $X_t-M_t$ conditioned on $\{M_t>\sqrt{2}t\}$.
\end{remark}

\section{Proof of Lemma \ref{initial-cond-V}}
In this section, we will give an upper estimate for $-\log\P_{\delta_x}(X_s([-A,A]^c=0, 0\le s\le t)$, which implies Lemma \ref{initial-cond-V}. Pinsky \cite{Pinsky} has proved a similar result
for super-Brownian motions with quadratic branching mechanism.
Here we use the  idea of \cite{Pinsky} to generalize the result to
super-Brownian motions with more general branching mechanisms.
\begin{lemma}{\bf (Maximum principle)}\label{l:eMP}
Let $\tl{\psi}(\lambda):=-a\lambda+b\lambda^{1+\vartheta}$, where $a>0, b>0, \vartheta>0$.
Assume that $v_1(x)$ and $v_2(x)$ are two functions defined on
$(a_1,a_2)$
such that $v_i(x)\ge (ab^{-1})^{1/\vartheta}$, $i=1,2$,
$v_1(a_i)\le v_2(a_i)$, $i=1,2$,
and that
$$
\frac{1}{2}\frac{d^2}{d x^2}v_{2}(x)-\tl\psi(v_2(x))\le \frac{1}{2}\frac{d^2}{dx^2}v_{1}(x)-\tl\psi(v_1(x)), \quad x\in(a_1,a_2).
$$
Then we have that
$$v_1(x)\le v_2(x),\quad x\in(a_1,a_2).$$
\end{lemma}
{\bf Proof:}
The proof is a slight modification of the proof of \cite[Proposition 3.1]{Bramson}, using \cite[Theorem 3.4]{PW}. See also the proof of \cite[Proposition 6.4]{Bovier}. We omit the details
\hfill$\Box$

\begin{lemma}\label{lem-h}
Let $\tl{\psi}(\lambda):=-a\lambda+b\lambda^{1+\vartheta}$, where $a>0, b>0, \vartheta\in(0,1]$.
For any $A>0$,
there exists an even function
$h_A(x)$ on $(-A,A)$ such that
\begin{align}\label{ODE-h}
 &\frac{1}{2}\Delta h_A(x)=\tl{\psi}(h_A(x)),\quad |x|<A,
\end{align}
and that $\lim_{x\to A}h_A(x)=\lim_{x\to -A}h_A(x)=\infty$.
Moreover,
there exist positive constants $c_1=c_1(a,b,\vartheta)$, $c_2=c_2(a,b,\vartheta)$ and $c_3=c_3(a,b,\vartheta)$ such that
 \begin{itemize}
   \item [(1)]$\max\{(ab^{-1})^{1/\vartheta},c_2A^{2/\vartheta}(A^2-x^2)^{-2/\vartheta}\}\le h_A(x)\le (ab^{-1})^{1/\vartheta}(1+c_1A^{2/\vartheta}(A^2-x^2)^{-2/\vartheta})$ for $|x|<A;$
   \item [(2)]$\frac{|h_A'(x)|}{h_A(x)}\le \frac{c_3}{A-|x|}$, for $|x|<A.$
 \end{itemize}
\end{lemma}
{\bf Proof:}
{\bf Step 1:}
First, for any
$m>(ab^{-1})^{1/\vartheta}$,
let $h_m(x)$ be the solution to the problem:
\begin{align}\label{ODE-hm}
 &\frac{1}{2}\Delta h_m(x)=\tl{\psi}(h_m(x)),\quad |x|<A,\\
 &h_m(A)=h_m(-A)=m.
\end{align}
Clearly $h_m$ is even. Since
$(ab^{-1})^{1/\vartheta}$ is a solution of $-a\lambda+b\lambda^{1+\vartheta}=0$,
the maximum principle in Lemma \ref{l:eMP} implies that
$h_m(x)\ge (ab^{-1})^{1/\vartheta}$ for $|x|< A$.

{\bf Step 2}
We want to find  $c_1>0$ such that the function $g(x)=(ab^{-1})^{1/\vartheta}(1+c_1A^{2/\vartheta}(A^2-x^2)^{-2/\vartheta})$ satisfies
$$
\frac{1}{2}\Delta g(x)\le \tl{\psi}(g(x))=-ag(x)+bg(x)^{1+\vartheta},\quad |x|<A.
$$
Assuming  this claim holds,  then using the  maximum principle
 in Lemma \ref{l:eMP} and the fact $\lim_{x\to A}g(x)=\lim_{x\to -A}g(x)=\infty$, we would get
$$g(x)\ge h_m(x), \quad |x|<A.$$
Now we prove the claim. Since
$\lim_{\lambda\downarrow0}
\frac{-(1+\lambda)+(1+\lambda)^{1+\vartheta}}{\lambda^{1+\vartheta}}=\infty$
and
$\lim_{\lambda\to\infty}
\frac{-(1+\lambda)+(1+\lambda)^{1+\vartheta}}{\lambda^{1+\vartheta}}=1$,
we have
$$
c_4:=
\inf_{\lambda\ge 0}
\frac{-(1+\lambda)+(1+\lambda)^{1+\vartheta}}{\lambda^{1+\vartheta}}
\in(0,\infty).
$$
Thus, we have that
\begin{align*}
-ag(x)+bg(x)^{1+\vartheta}&\ge c_4a(ab^{-1})^{1/\vartheta}c_1^{1+\vartheta}A^{2+2/\vartheta}(A^2-x^2)^{-2(1+\vartheta)/\vartheta}.
\end{align*}
It is clear that, for any $x\in(-A,A)$,
\begin{align*}
\frac{1}{2}\Delta g(x)&=(ab^{-1})^{1/\vartheta}c_1A^{2/\vartheta}2\vartheta^{-1}[A^2+(4\vartheta^{-1}+1)x^2](A^2-x^2)^{-2-2/\vartheta}\\
&\le (ab^{-1})^{1/\vartheta}c_12\vartheta^{-1}(4\vartheta^{-1}+2)A^{2+2/\vartheta}(A^2-x^2)^{-2-2/\vartheta}.
\end{align*}
Therefore it suffices to choose
$$
c_1=\Big(4c_4^{-1}a^{-1}\vartheta^{-1}(\frac{2}{\vartheta}+1)\Big)^{1/\vartheta}.$$

{\bf Step 3}
For any $\delta>0$, define $g_\delta(x):=\frac{c_2A^{2/\vartheta}}{((A+\delta)^2-x^2)^{2/\vartheta}}$, where $c_2>(ab^{-1})^{1/\vartheta}$ is a constant. We claim that there exists $c_2=c_2(a,b,\vartheta)>0$ such that
\begin{equation}\label{4.9}
 \frac{1}{2}\Delta g_\delta(x)\ge -ag_\delta(x)+bg_\delta(x)^{1+\vartheta},\quad |x|<A+\delta.
\end{equation}
Assuming  this claim holds, then
applying the  maximum principle in Lemma \ref{l:eMP}, we would get that, for $m$ large enough
$$h_m(x)\ge g_\delta(x),\quad |x|<A.$$
Now we prove the claim. In fact,
$$
\frac{1}{2}\Delta g_\delta(x)\ge c_22\vartheta^{-1}A^{2+2/\vartheta}((A+\delta)^2-x^2)^{-2-2/\vartheta}
$$
and
$$
-ag_\delta(x)+bg_\delta(x)^{1+\vartheta}\le bg_\delta(x)^{1+\vartheta}=bc_2^{1+\vartheta}A^{2/\vartheta+2}((A+\delta)^2-x^2)^{-2-2/\vartheta}.
$$
Thus we only need to choose
$$
c_2=(2b^{-1}\vartheta^{-1})^{1/\vartheta}.
$$

{\bf Step 4}
By the maximum principle in Lemma \ref{l:eMP},  $h_m$ is non-decreasing in $m$, thus $h_A(x):=\lim_{m\to\infty} h_m(x)$ exists. Hence for any $\delta>0$,
$$g_\delta(x)\le h_A(x)\le g(x).$$
Letting $\delta\to 0$, we have that, for any $|x|<A$,
$$
\frac{c_2A^{2/\vartheta}}{(A^2-x^2)^{2/\vartheta}}\le h_A(x)\le (ab^{-1})^{1/\vartheta}(1+c_1A^{2/\vartheta}(A^2-x^2)^{-2/\vartheta}).
$$
Clearly $\lim_{x\to A}h_A(x)=\lim_{x\to -A}h_A(x)=\infty$.

{\bf Step 5}
Now we show that $h_A$ satisfies \eqref{ODE-h}.
By \eqref{ODE-hm}, we have that for any $0<A'<A$,
$$h_m(x)=
    -\mE_x\int_0^{\tau_{A'}}\tl\psi(h_m(B_s))\,ds+\mE_x(h_m(B_{\tau_{A'}})),\quad x\in(-A',A'),
$$
where $\tau_{A'}$ is the exit time of $B$ from $(-A',A')$.
Letting  $m\to\infty$ and applying the dominated convergence theorem,  we get that
$$
h_A(x)=
    -\mE_x\int_0^{\tau_{A'}}\tl\psi(h(B_s))\,ds+\mE_x(h(B_{\tau_{A'}})),\quad x\in(-A',A'),
$$
which implies that $h_A$ satisfies \eqref{ODE-h} for $x\in(-A',A')$.
Since $A'\in (0, A)$ is arbitrary,
$h_A$ satisfies \eqref{ODE-h} for $x\in(-A,A)$.

{\bf Step 6} Finally, we prove that
$\frac{|h_A'(x)|}{h_A(x)}\le \frac{c_3}{A-|x|},\quad |x|<A.$
Since $h_A$ is an even function, we have $\frac{|h_A'(x)|}{h_A(x)}=\frac{|h_A'(|x|)|}{h_A(|x|)}$. To prove the desired result, we only need to consider $x\ge 0$.
Since $h_A(x)\ge (ab^{-1})^{1/\vartheta}$ and
$$
\frac{1}{2}\Delta h_A(x)=\tl{\psi}(h_A(x))\ge 0, \quad |x|<A,
$$
we know that $h_A'(x)$ is increasing on $(-A,A)$.
Since  $h_A$ is an even function, we have $h_A'(0)=0$. Thus, $h_A'(x)\ge0$, for $x\in[0,A)$£¬ which implies that
\begin{equation}\label{4.12}
  \frac{h_A'(x)}{h_A(x)}\ge0,\quad x\in[0,A).
\end{equation}
Define $w_1(x)=\frac{2a(c_1)^\vartheta}{A-x}-\frac{h_A'(x)}{h_A(x)}$, for $x\in[0,A).$ Then, for any $x\in(0,A)$,
\begin{align*}
  w_1'(x)=\frac{2(c_1)^\vartheta}{(A-x)^2}-2(bh_A(x)^\vartheta-a)+\Big(\frac{h_A'(x)}{h_A(x)}\Big)^2\ge 0,
\end{align*}
where the last inequality follows from the fact that
  $$
  bh_A(x)^\vartheta-a\le a(1+c_1A^{2/\vartheta}(A^2-x^2)^{-2/\vartheta})^\vartheta-a\le ac_1^{\vartheta}A^{2}(A^2-x^2)^{-2}\le \frac{a(c_1)^\vartheta}{(A-x)^2}.
  $$
Since $h_A'(0)=0$, then $w_1(0)>0$.  Thus for any $x\in(0,A)$, $w_1(x)\ge w_1(0)>0$, that is
\begin{equation}\label{5.2}
  \frac{2a(c_1)^\epsilon}{A-x}\ge\frac{h_A'(x)}{h_A(x)},\quad x\in[0,A).
\end{equation}
Combining \eqref{4.12} and \eqref{5.2}, we get the desired result.
\hfill$\Box$

\begin{lemma}\label{lem:com2}
Assume ${\bf (H1)}$ and ${\bf (H3)}$ hold.
Then, there exist positive constants
$c_4=c_4(a,b,\vartheta)$ and  $c_5=c_5(a,b,\vartheta)$ such that  for any $A>0$ and $|x|<A$,
\begin{equation}\label{5.3}
  -\log\P_{\delta_x}(X_s([-A,A]^c)=0,\forall s\in[0,t])\le h_A(x)\exp\left\{-\left(c_4\frac{(A-|x|)^2}{t}-at-c_5\right)\right\}.
\end{equation}
\end{lemma}
{\bf Proof:}
Let $\tilde{X}$ be a super-Brownian motion with branching mechanism $\tl{\psi}(\lambda):=-a\lambda+b\lambda^{1+\vartheta}$.
Define
$$
h(t,x):=-\log\P_{\delta_x}(\tl{X}_s([-A,A]^c)=0,\forall s\le t)
$$
and
$$
h_m(t,x):=-\log\E_{\delta_x}\Big[\exp\Big\{-\int_0^t \langle \phi_m,\tl X_s\rangle\,ds\Big\}\Big],
$$
where $\phi_{m}\in C^\infty(\R)$ satisfies
\begin{align*}
         \phi_m(y)=0, & \quad |y|<A, |y|>A+m+1, \\
         \phi_m(y)=m, & \quad A+\frac{1}{m}\le |y|\le A+m.
\end{align*}
Then
$h(t,x)=\lim_{m\to\infty}h_m(t,x)$
and $h_m(t,x)$ satisfies the equation
$$h_m(t,x)+\mE_x\int_0^t\tl\psi(h_m(t-s,B_s))\,ds=\mE_x\int_0^t \phi_m(B_s)\,ds.$$
For the display above we refer the readers to \cite[Corollay 5.17]{Li11}. Thus,
\begin{align*}
 &\frac{\partial h_m}{\partial t}(t,x)-\frac{1}{2}\Delta h_m(t,x)=-\tl{\psi}(h(t,x))+\phi_m(x),\quad t>0,
\end{align*}
which implies that
\begin{align*}
 &\frac{\partial h_m}{\partial t}(t,x)-\frac{1}{2}\Delta h_m(t,x)=-\tl{\psi}(h(t,x)),\quad |x|<A, t>0.
\end{align*}
Since $\psi\ge \tl{\psi}$, then using  arguments similar to  that used in \cite[Corollary 5.18]{Li11}, we get that
 $$-\log\E_{\delta_x}\Big[\exp\Big\{-\int_0^t \langle \phi_m, X_s\rangle\,ds\Big\}\Big]\le h_m(t,x).$$
Letting $m\to\infty$, we get
$$-\log\P_{\delta_x}(X_s([-A,A]^c)=0,\forall s\in[0,t])\le h(t,x),$$
so it suffices to show that the result holds for $h(t,x)$.

Let $f$ be an even function satisfying
\begin{align}\label{4.7}
  &f\in C^2([-1,1]),\quad f(y)>0, \mbox{ if } -1< y<1;\nonumber\\
  & f(0)=1,\quad f'(0)=0,\quad f(1)=0,\quad f'(1)=0,\quad f''(1)>0.\nonumber\\
  &\sup_{y\in[0,1]}\frac{(f'(y))^2}{f(y)}<\infty.
\end{align}
It has been proved in the proof of  \cite[Theorem 1]{Pinsky} that such $f$ exists.
Define
$$
v(t,x):=h_A(x)\exp\Big\{c_5+at-\frac{\delta A^2}{t}f(\frac{x}{A}))\Big\},
\quad |x|<A,
$$
where $c_5,\delta>0$ are to be fixed later.
It is clear that $\lim_{t\to 0}v(t,x)=0,$ $\lim_{|x|\to A}v(t,x)=\infty,$
since $\lim_{|x|\to A}h_A(x)=\infty.$

To prove the result, we want to find suitable
$c_{5},\delta$ such that
\begin{equation}\label{4.8}
  \frac{\partial v}{\partial t}(t,x)-\frac{1}{2}\Delta v(t,x)\ge av(t,x) -bv(t,x)^{1+\vartheta},\quad |x|<A.
\end{equation}
Assuming this claim for the time being, by the maximum principle
 in Lemma \ref{l:mp}, we would have $h_m(t,x)\le v(t,x).$ Letting $m\to\infty$, we get
$$
h(t,x)\le v(t,x),\quad |x|<A,t>0.
$$
Since $f''(1)>0$, we have $\inf_{y\in[0,1]}\frac{f(y)}{(1-y)^2}>0$. Thus,
$$
h(t,x)\le v(t,x)\le h_A(x)\exp\Big\{c_5+at-\frac{c_4 (A-|x|)^2 }{t}\Big\},
$$
where $c_4=\delta\inf_{y\in[0,1]}\frac{f(y)}{(1-y)^2}>0$.

Now we prove \eqref{4.8}.
Note that, by \eqref{ODE-h}, \eqref{4.8} is equivalent to, for $x\in[0,A)$, $y=x/A$,
$$
a+\frac{\delta A^{2}}{t^2}f(y)-\frac{\delta^2 A^{2}}{2t^2}(f'(y))^2+\frac{h_A'(x)}{h_A(x)}\frac{\delta A}{t}f'(y)+\frac{\delta }{2t}f''(y)\ge -bv(t,x)^{\vartheta}+bh_A(x)^{\vartheta}.
$$
Note that $\frac{(f'(y))^2}{f(y)}$, $\frac{|f'(y)|}{1-y}$, and $f''(y)$ are all bounded. Let $K$ be the common upper bound. By Lemma \ref{lem-h}, $\frac{|h_A'(x)|}{h_A(x)}\le c_3(A-x)^{-1}$.
Choose $\delta\in (0, K^{-1})$.
It suffices to show that
\begin{equation}\label{4.11}
  a+\frac{\delta A^{2}}{2t^2}f(y)-\frac{c_3K\delta}{t}-\frac{K\delta}{2t}\ge -bv(t,x)^{\vartheta}+bh_A(x)^{\vartheta}.
\end{equation}
If $\frac{\delta A^{2}}{t}f(y)\ge c_5/2$, then
the left hand side of  \eqref{4.11} is bigger than
$$
a+\frac{c_5}{4t}-\frac{c_3}{t}-\frac{1 }{2t},$$
and by Lemma \ref{lem-h}, the right hand side of \eqref{4.11} is less than
\begin{align*}
bh_A(x)^{\vartheta}&\le a(1+c_1A^{2/\vartheta}(A^2-x^2)^{-2/\vartheta})^{\vartheta}\le a(1+c_1^{\vartheta}A^{-2}(1-y^2)^{-2})\\
&=a+ac_1^{\vartheta}(f(y))^{-1}A^{-2}\frac{f(y)}{(1-y^2)^{2}}\le a+a\frac{2\delta c_1^{\vartheta}K}{c_5 t}\le a+\frac{2ac_1^{\vartheta}}{c_5 t}.
\end{align*}
Thus, when we choose $c_5$ large enough, \eqref{4.11} is true.

If $c_5/2\ge \frac{\delta A^{2}}{t}f(y)$, then
the left hand side of  \eqref{4.11} is bigger than
$$
a-\frac{c_3}{t}-\frac{1}{2t},
$$
and the right hand side of \eqref{4.11} is less than
\begin{align*}
bh_A(x)^{\vartheta}(1-e^{\vartheta c_5/2})&\le-(e^{\vartheta c_5/2}-1)bc_2^{\vartheta}A^2(A^2-x^2)^{-2}\\
&= -bc_2^{\vartheta}(e^{\vartheta c_5/2}-1)\frac{1}{f(y)A^{2}}\frac{f(y)}{(1-y)^2(1+y)^2}\\
&\le -\frac{bc_2^{\vartheta}\delta(e^{\vartheta c_5/2}-1)}{2 c_5}\inf_{y\in[0,1]}\frac{f(y)}{(1-y)^2} \frac{1}{t}.
\end{align*}
Since $\inf_{y\in[0,1]}\frac{f(y)}{(1-y)^2}>0$, we can choose $c_{5}$ large enough such that \eqref{4.11} is true.
The proof is now complete.
\hfill $\Box$

\noindent{\bf Proof of Lemma \ref{initial-cond-V}:}
It is clear that
$$\E\Big[\exp\Big\{-\int_\R \phi(y-x)X_t(dy)\Big\} , M_t\le x\Big]\ge \P(\|X_t\|=0)>0,$$
where the last inequality follows from \eqref{extinction}. Thus
$$V(t,x)\le -\log\P(\|X_t\|=0)<\infty,$$ which implies that
$V(t,\cdot)$ is a bounded function.
For any $x>1$, it follows from  Lemma \ref{lem:com2} that
\begin{align*}
  V(t,x) & \le -\log \P(X_s([-x,x]^c)=0,s\le t)\le
   h_x(0)\exp\{-\frac{c_4}{t}x^2+at+c_5\} \\
  &\le c(t) e^{-c_4x^2/t},
\end{align*}
where $c(t)$ is a constant which may depend on $t$.
Thus, the desired result follows.
\hfill
$\Box$

\bigskip
\noindent
{\bf Acknowledgment:}
	We thank the  referee for helpful comments.

\vspace{.1in}
\begin{singlespace}

\end{singlespace}

\end{doublespace}

\vskip 0.2truein
\vskip 0.2truein

\noindent{\bf Yan-Xia Ren:} LMAM School of Mathematical Sciences \& Center for
Statistical Science, Peking
University,  Beijing, 100871, P.R. China. Email: {\texttt
yxren@math.pku.edu.cn}

\smallskip
\noindent {\bf Renming Song:} Department of Mathematics,
University of Illinois,
Urbana, IL 61801, U.S.A.
Email: {\texttt rsong@illinois.edu}

\smallskip

\noindent{\bf Rui Zhang:}  School of Mathematical Sciences \& Academy for Multidisciplinary Studies, Capital Normal
University,  Beijing, 100048, P.R. China. Email: {\texttt
zhangrui27@cnu.edu.cn}

\end{document}